\documentstyle[bezier]{article}

\makeatletter
  
  \@addtoreset{equation}{section}
\makeatother

\newtheorem{theorem}{Theorem}[section]
\newtheorem{proposition}{Proposition}[section]
\newtheorem{remark}{Remark}[section]
\newtheorem{example}{Example}[section]
\newtheorem{conjecture}{Conjecture}[section]

\newfont{\germ}{eufm10}
\newfont{\slsmall}{cmsl8}

\def\B{{\cal B}}
\def\bbar{\overline{b}}
\def\cbox#1#2{\fbox{$\hspace{-.8mm}{#1\atop#2}\hspace{-.8mm}$}}
\def\cboxfig#1#2{\framebox(10,20){$\mbox{#1}\atop\mbox{#2}$}}
\def\cd{\cdots}
\def\eps{\epsilon}
\def\et#1{\tilde{e}_{#1}}
\def\ft#1{\tilde{f}_{#1}}
\def\geh{\goth{g}}
\def\goth#1{\mbox{\germ #1}}
\def\L{{\cal L}}
\def\La{\Lambda}
\def\la{\lambda}
\def\ol#1{\overline{#1}}
\def\ot{\otimes}
\def\P{{\cal P}}
\def\pbar{\overline{p}}
\def\Pcl{P_{cl}}
\def\Pcll{(P_{cl}^+)_l}
\def\Proof{\noindent{\sl Proof.}\quad}
\def\Q{{\bf Q}}
\def\qed{~\rule{1mm}{2.5mm}}
\def\Uq{U_q(\geh)}
\def\veps{\varepsilon}
\def\vphi{\varphi}
\def\wt{\mbox{\sl wt}\,}
\def\Z{{\bf Z}}

\begin{document}

\title{ Crystals for Demazure Modules of \\
        Classical Affine Lie Algebras }

\author{
Atsuo Kuniba\thanks{
Institute of Physics, University of Tokyo, Komaba, Tokyo 153, Japan}, 
Kailash C. Misra\thanks{
Department of Mathematics, North Carolina State University, 
Raleigh, NC 27695-8205, USA}, 
Masato Okado\thanks{
Department of Mathematical Science, Faculty of Engineering Science,
Osaka University, Toyonaka, Osaka 560, Japan},\\
Taichiro Takagi\thanks{
Department of Mathematics and Physics, National Defense Academy,
Yokosuka 239, Japan}
and Jun Uchiyama\thanks{
Department of Physics, Rikkyo University, Nishi-Ikebukuro, Tokyo 171, Japan}
} 

\date{}
\maketitle

\begin{abstract}
\noindent
We study, in the path realization, crystals for Demazure modules
of affine Lie algebras of types $A^{(1)}_n,B^{(1)}_n,C^{(1)}_n,D^{(1)}_n,
A^{(2)}_{2n-1},A^{(2)}_{2n}, and D^{(2)}_{n+1}$. We find a 
special sequence of affine Weyl group elements for the selected 
perfect crystal, and show if the highest weight is $l\La_0$, the
Demazure crystal has a remarkably simple structure.
\end{abstract} 

\setcounter{section}{-1}

\section{Introduction}

Let $\geh$ be a symmetrizable Kac-Moody algebra. Let $\Uq$ be 
its quantized universal enveloping algebra and $V(\la)$ be the 
integrable $\Uq$-module with dominant integral highest weight
$\la$. Let $W$ be the Weyl group of $\geh$. Fixing an element 
$w$ of $W$, the Demazure module $V_w(\la)$ is defined as a
finite dimensional subspace of $V(\la)$ generated from the 
extremal weight space $V(\la)_{w\la}$ by the $e_i$ generators of 
$\Uq$. Though the Demazure module itself can be defined in the 
same way for the classical case $q=1$, we stay in the quantum
case. The reason is the existence of a ``good" basis. Let $(\L(\la),
\B(\la))$ be the crystal base of $V(\la)$ \cite{K1}. In \cite{K2},
Kashiwara showed there exists a subset $\B_w(\la)$ of $\B(\la)$
such that 
\[
V_w(\la)=\bigoplus_{b\in\B_w(\la)}\Q(q)G_\la(b).
\]
Here $G_\la(b)$ is the lower global base \cite{K3}. Moreover,
$\B_w(\la)$ has a quite remarkable simple recursive property:
\[
\mbox{If }r_iw\succ w,
\mbox{ then }\B_{r_iw}(\la)
=\bigcup_{n\ge0}\ft{i}^n\B_w(\la)\setminus\{0\}.
\]
Here $r_i$ is a simple reflection of $W$ and $\succ$ denotes the 
Bruhat order. In fact, using this property Kashiwara gave a new proof
of Demazure's character formula for an arbitrary symmetrizable 
Kac-Moody algebra \cite{K2}.

Let us now focus on the quantum affine algebra $\Uq$, where $\geh$
is of affine type. In this case, we have a description of $\B(\la)$ in
terms of paths \cite{KMN1,KMN2}. Roughly speaking, the set of paths
is a suitable subset of the half infinite tensor product of a ``perfect"
crystal $B$ which is a crystal of a finite dimensional $U'_q(\geh)$-module 
having some nice properties \cite{KMN1}. On this set, the actions of
$\et{i}$ and $\ft{i}$ are explicitly given. In \cite{KMOU}, we gave a 
criterion for $\B_w(\la)$ to have a tensor product structure. To
describe the general situation, the mixing index $\kappa$ was introduced. 
Taking $\kappa=1$ for simplicity, the result in \cite{KMOU} is stated
as follows. Consider an increasing sequence $\{w^{(k)}\}_{k\ge0}$
of $W$ with respect to the Bruhat order. If a perfect crystal $B$ and
$\{w^{(k)}\}$ satisfy several assumptions, then $\B_{w^{(k)}}(\la)$ is
given by
\begin{equation} \label{eq:intro}
\begin{array}{ccccccccc}
\B(\la)&\subset&\cd\ot&B&\ot&B&\ot&B&\ot B\ot\cd\ot B\\
\cup&&&&&&&&\\
\B_{w^{(k)}}(\la)&=&
\cd\ot&\hspace{-2.5mm}\bbar_{j+2}\hspace{-2.5mm}
&\ot&\hspace{-2.5mm}\bbar_{j+1}\hspace{-2.5mm}
&\ot&\hspace{-2.5mm}B^{(j)}_a\hspace{-2.5mm}&\ot B\ot\cd\ot B.
\end{array}
\end{equation}
Here $j,a$ are determined from $k$, $B^{(j)}_a$ is a subset of $B$, and
$\bbar=\cd\ot\bbar_j\ot\cd\ot\bbar_1$ is the ground state path corresponding 
to the highest weight vector in $\B(\la)$.

The purpose of this article is to show that if $\la=l\La_0$ (we also
discuss some other similar cases), we can find the sequence $\{w^{(k)}\}$
satisfying (\ref{eq:intro}) for $\geh$ of classical types 
(i.e. $\geh=A^{(1)}_n,B^{(1)}_n,D^{(1)}_n,A^{(2)}_{2n-1},A^{(2)}_{2n},
D^{(2)}_{n+1},C^{(1)}_n$).
We choose the perfect crystal $B$ from the list in \cite{KMN2} except
for the $C^{(1)}_n$ case. To illustrate, we take an example of 
$\geh=A^{(1)}_3,B=B(\La_2)$. The crystal graph of $B$ is given
in Figure 1.
\begin{figure}[ht]
\vskip1cm
\begin{picture}(300,130)(0,0)
\put(50,0){
\put(5,65){\cboxfig12}
\put(20,75){\vector(1,0){40}}
\put(34,80){\makebox(10,10){\bf 2}}
\put(65,65){\cboxfig13}
\put(115,105){\cboxfig23}
\put(80,79){\vector(1,1){31}}
\put(87,97){\makebox(10,10){\bf 1}}
\put(80,71){\vector(1,-1){31}}
\put(93,57){\makebox(10,10){\bf 3}}
\put(115,25){\cboxfig14}
\put(130,111){\vector(1,-1){31}}
\put(143,97){\makebox(10,10){\bf 3}}
\put(130,39){\vector(1,1){31}}
\put(139,57){\makebox(10,10){\bf 1}}
\put(165,65){\cboxfig24}
\put(180,75){\vector(1,0){40}}
\put(194,80){\makebox(10,10){\bf 2}}
\put(225,65){\cboxfig34}
\bezier{500}(170,58)(90,-50)(12,58)
\put(12,58){\vector(-1,2){1}}
\bezier{500}(230,92)(140,200)(72,92)
\put(72,92){\vector(-1,-2){1}}
\put(85,9){\makebox(10,10){\bf 0}}
\put(145,148){\makebox(10,10){\bf 0}}
}
\end{picture}
\caption{Level 1 perfect crystal $B(\Lambda_2)$ for $A^{(1)}_3$}
\end{figure}
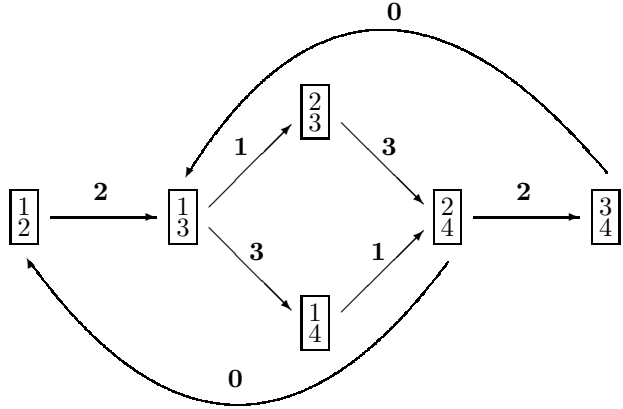
Since $B$ is a level 1 perfect crystal,
let us take $\la=\La_0$. Then the path $\bbar$ is given by 
\[
\bbar=\cd\ot\cbox12\ot\cbox34\ot\cbox12\ot\cbox34.
\]
In this case, the sequence of Weyl group elements $\{w^{(k)}\}_{k\ge0}$
is given as follows:
\[
w^{(0)}=1,\qquad w^{(k+1)}=r_iw^{(k)}\quad(k\ge0),
\]
where $i=0$ ($k\equiv0,3$), $i=3$ ($k\equiv1,6$), $i=1$ ($k\equiv2,5$), 
$i=2$ ($k\equiv4,7$). Here $\equiv$ denotes `congruence modulo 8'. The
integers $j,a$ in (\ref{eq:intro}) are determined from $k$ by 
$k=4(j-1)+a,j\ge1,0\le a<4$, and $B^{(j)}_a$ are given as follows:

\noindent If $j$ is odd,
\[
B^{(j)}_0=\left\{\cbox34\right\},\;\;
B^{(j)}_1=B^{(j)}_0\cup\left\{\cbox13\right\},\;\;
B^{(j)}_2=B^{(j)}_1\cup\left\{\cbox14\right\},\;\;
B^{(j)}_3=B^{(j)}_2\cup\left\{\cbox23,\cbox24\right\}.
\]

\noindent If $j$ is even,
\[
B^{(j)}_0=\left\{\cbox12\right\},\;\;
B^{(j)}_1=B^{(j)}_0\cup\left\{\cbox13\right\},\;\;
B^{(j)}_2=B^{(j)}_1\cup\left\{\cbox23\right\},\;\;
B^{(j)}_3=B^{(j)}_2\cup\left\{\cbox14,\cbox24\right\}.
\]

This paper is organized as follows. In the next section, we review
perfect crystals, Demazure crystals and the criterion developed in 
\cite{KMOU}. From section 2 to 8, we list the sequence of Weyl group
elements $\{w^{(k)}\}$, the subset $B^{(j)}_a\subset B$ and other
important data for Demazure crystals. In the last section, some 
observations and conjectures are given.

\section{Crystals for Demazure modules}

In \cite{KMOU}, a criterion for Demazure crystals is presented.
It clarifies their tensor product structure, and involves a parameter
$\kappa$ which measures the degree of {\em mixing}. Here we summarize
the criterion for $\kappa=1$.

\subsection{Perfect crystal}
We follow the notations of the quantized universal enveloping algebra
and the crystal base in \cite{KMN1}, except that we shall use a different
font for the crystal base $(\L(\la),\B(\la))$ of the irreducible highest
weight module $V(\la)$ in order to avoid the confusion with the crystal
base $(L,B)$ of $V$ in $\mbox{Mod}^f(\geh,\Pcl)$ which may also have an
argument. We review necessary properties of perfect crystals. Our main
reference is \cite{KMN1}.

Let $B$ be a perfect crystal of level $l$. Then for any $\la\in\Pcll$, 
there exists a unique element $b(\la)\in B$ such that $\vphi(b(\la))=\la$.
Let $\sigma$ be the automorphism of $\Pcll$ given by $\sigma\la=\veps(b(\la))$.
We set $\bbar_k=b(\sigma^{k-1}\la)$ and $\la_k=\sigma^k\la$. Then perfectness
assures the following isomorphism of crystals.
\begin{equation} \label{eq:iteration}
\B(\la_{k-1})\simeq\B(\la_k)\ot B.
\end{equation}
Define the set of paths $\P(\la,B)$ by
\[
\P(\la,B)=\{p=\cd\ot p(2)\ot p(1)\mid p(j)\in B,p(k)=\bbar_k\mbox{ for }k\gg1\}.
\]
Iterating the isomorphism (\ref{eq:iteration}), we see $\B(\la)$ is isomorphic 
to $\P(\la,B)$. Under this isomorphism, the highest weight vector $u_\la$ 
in $\B(\la)$ corresponds to the path $\pbar=\cd\ot\bbar_k\ot\cd\ot\bbar_2\ot
\bbar_1$, which we call the {\em ground state} path. The actions of $\et{i}$
and $\ft{i}$ on $\P(\la,B)$ are determined by {\em signature rule}, which
we explain in the next subsection.

\subsection{Signature rule} \label{subsec:signature}
We need to know the actions of $\et{i}$ and $\ft{i}$ on the set of paths 
$\P(\la,B)$. To see this, we consider the following isomorphism. 
\begin{equation}\label{eq:iso1}
\P(\la,B)\simeq \B(\la_k)\ot B^{\ot k}.
\end{equation}
For each $p\in\P(\la,B)$, if we take $k$ sufficiently large, we can assume
that $p$ corresponds to $u_{\la_k}\ot p(k)\ot\cdots\ot p(1)$ with $u_{\la_k}$
the highest weight vector of $\B(\la_k)$ and $p(j)\in B$ ($j=1,\cdots,k$). 
Then we apply Proposition 2.1.1 in \cite{KN} to see on which component $\et{i}$
or $\ft{i}$ acts. Let us suppose that $\et{i}$ or $\ft{i}$ acts on the $j$-th
component from the right end. If $j<k+1$, we have
\begin{equation}\label{eq:action_e}
\et{i}p=\cdots\ot\et{i}p(j)\ot\cdots\ot p(2)\ot p(1)
\end{equation}
or
\begin{equation}\label{eq:action_f}
\ft{i}p=\cdots\ot\ft{i}p(j)\ot\cdots\ot p(2)\ot p(1).
\end{equation}
If $j=k+1$, we see we should have taken $k$ larger. This happens only for 
$\ft{i}$.

This determination of the component can be rephrased using the notion of
{\em signature}. Let $p\in\P(\la,B)$ correspond to $u_{\la_k}\ot p(k)\ot
\cdots\ot p(1)$ under (\ref{eq:iso1}) as above. With $p(j)$ ($1\le j\le k$), 
we associate 
\begin{eqnarray*}
\eps^{(j)}&=&(\eps^{(j)}_1,\cdots,\eps^{(j)}_m),\\
m&=&\veps_i(p(j))+\vphi_i(p(j)),\\
\eps^{(j)}_a&=&-,\:\:\mbox{if}\:\: 1\le a\le\veps_i(p(j)),
\quad +,\:\:\mbox{if}\:\: \veps_i(p(j))<a\le m.
\end{eqnarray*}
For the highest weight vector $u_{\la_k}$, we take 
\[
\eps^{(k+1)}=(\underbrace{+,\cdots,+}_{\langle\la_k,h_i\rangle}).
\]
We then append these $\eps^{(j)}$'s so that
we have
\[
\eps=(\eps^{(k+1)},\eps^{(k)},\cdots,\eps^{(1)}).
\]
We call it ($i$-)signature of $p$ truncated at the $k$-th position.

Next we consider a sequence of signatures.
\[
\eps=\eta_0,\eta_1,\cdots,\eta_{\max}.
\]
Here $\eta_{j+1}$ is obtained from $\eta_j$ by deleting the leftmost
adjacent $(+,-)$ pair of $\eta_j$. Eventually, we arrive at the 
following signature
\[
\eta_{\max}=(\underbrace{-,\cdots,-}_{n_-},\underbrace{+,\cdots,+}_{n_+}),
\]
with $n_\pm\ge0$. We call it the {\em reduced signature} and denote by
$\ol{\eps}$. 

The component on which $\et{i}$ or $\ft{i}$ acts in (\ref{eq:action_e})
or (\ref{eq:action_f})
reads as follows. If $n_-=0$ (resp. $n_+=0$), we set $\et{i}p=0$ 
(resp. $\ft{i}p=0$). Otherwise, take the rightmost $-$ (resp. leftmost $+$), 
and find the component $\eps^{(j)}$ to which it belonged. Then, this $j$
is the position in (\ref{eq:action_e}) (resp. (\ref{eq:action_f}))
we looked for. Note that if $k$ is large enough, the
position $j$ does not depend on the choice of $k$. 

\begin{remark}
Of course, this signature rule can be applied to the tensor product
of crystals $B_1\ot\cd\ot B_l$.
\end{remark}

\begin{example}
Let $\geh=A^{(1)}_1$, $B$ be the classical crystal of
the 3-dimensional irreducible representation. Its crystal graph is 
described as follows.
\[
B\hspace{5mm}00\mathop{\rightleftharpoons}_0^1
01\mathop{\rightleftharpoons}_0^1 11
\]
It is known that $B$ is perfect of level 2. We have isomorphisms
$B(\la)\simeq\P(\la,B)$ for $\la=2\La_0,\La_0+\La_1,2\La_1$.
Let $\la=2\La_0$. We see $\la_k=2\La_0$ $(k:\mbox{even})$, $2\La_1$
$(k:\mbox{odd})$. The ground-state path is given by
\[
\ol{p}=\cdots\ot11\ot00\ot11\ot00\ot11.
\]
Consider a path
\[
p=\cdots\ot11\ot01\ot01\ot01\ot00.
\]
The dotted part of $p$ is the same as that of $\ol{p}$. Then the
1-signature of $p$ truncated at the 5-th position and its reduced 
signature read as follows.
\begin{eqnarray*}
\eps&=&(++,--,-+,-+,-+,++),\\
\ol{\eps}&=&(\mathop{-}^4\mathop{+}^2\mathop{+}^1\mathop{+}^1).
\end{eqnarray*}
Here the number above each sign shows the component to which it belonged.
Consequently, we have 
\begin{eqnarray*}
\et{1}p&=&\cdots\ot11\ot00\ot01\ot01\ot00,\\
\ft{1}p&=&\cdots\ot11\ot01\ot01\ot11\ot00.
\end{eqnarray*}
\end{example}

\subsection{Demazure crystal}
Let $\{r_i\}_{i\in I}$ be the set of simple reflections, and let
$W$ be the Weyl group. For $w\in W$, $l(w)$ denote the length of 
$w$, and $\prec$ denote the Bruhat order on $W$. Let $U_q^+(\geh)$
be the subalgebra of $\Uq$ generated by $e_i$'s. We consider the
irreducible highest weight $\Uq$-module $V(\la)$ ($\la\in\Pcl^+$).
Let $V_w(\la)$ denote the $U_q^+(\geh)$-module generated by the 
extremal weight space $V(\la)_{w\la}$. These modules
$V_w(\la)$ ($w\in W$) are called the Demazure modules.
Let $(\L(\la),\B(\la))$ be the crystal base of $V(\la)$. In \cite{K2}
Kashiwara showed that for each $w\in W$, there exists a subset
$\B_w(\la)$ of $\B(\la)$ such that
\[
\frac{V_w(\la)\cap\L(\la)}{V_w(\la)\cap q\L(\la)}
=\bigoplus_{b\in\B_w(\la)}\Q b.
\]
Furthermore, $\B_w(\la)$ has the following recursive property.
\begin{eqnarray}
&&\mbox{If }r_iw\succ w,\mbox{ then}\nonumber\\
&&\B_{r_iw}(\la)=\bigcup_{n\ge0}\ft{i}^n\B_w(\la)\setminus\{0\}. 
\label{recursive}
\end{eqnarray}
We call $\B_w(\la)$ the {\em Demazure crystal} associated with the 
Demazure module $V_w(\la)$.

\subsection{Criterion} \label{subsec:crit}
In this sebsection, we review the result of \cite{KMOU} when
$\kappa=1$. (For the definition of the {\em mixing index} $\kappa$,
see section 2.3 of \cite{KMOU}.)

Let $\la$ be an element of $\Pcll$, and let $B$ be a classical crystal. 
To state the result, we need to assume four conditions (I-IV).
\begin{description}
\item[(I)   ]
$B$ is perfect of level $l$.
\end{description}
Thus, we can assume an isomorphism between $\B(\la)$ and the set of paths
$\P(\la,B)$. Let $\pbar=\cd\ot\bbar_2\ot\bbar_1$ denote the ground state 
path. Fix a positive integer $d$. For a set of elements
$i_a^{(j)}$ ($j\ge1,1\le a\le d$) in $I$, we define
$B_a^{(j)}$ ($j\ge1,0\le a\le d$) by 
\[
B^{(j)}_0=\{\bbar_j\},\hspace{1cm}
B_a^{(j)}=\bigcup_{n\ge0}\ft{i_a^{(j)}}^n B_{a-1}^{(j)}\setminus\{0\}
\quad(a=1,\cdots,d).
\]
\begin{description}
\item[(II)  ]
For any $j\ge1$,
$B_d^{(j)}=B$.
\item[(III) ]
For any $j\ge1$ and $1\le a\le d$,
$\langle\la_j,h_{i^{(j)}_a}\rangle\le\veps_{i^{(j)}_a}(b)$
for all $b\in B^{(j)}_{a-1}$.
\end{description}
We now define an element $w^{(k)}$ of the Weyl group $W$ by
\[
w^{(0)}=1,\hspace{1cm}
w^{(k)}=r_{i^{(j)}_a}w^{(k-1)}\quad\mbox{for }k>0.
\]
Here, if $k=0$, we set $j=1,a=0$, and otherwise, $j$ and $a$ are 
fixed from $k$ by $k=(j-1)d+a,j\ge1,1\le a\le d$.
\begin{description}
\item[(IV)  ]
$w^{(0)}\prec w^{(1)}\prec\cd\prec w^{(k)}\prec\cd$.
\end{description}
We shall discuss how to check this condition in the next subsection.

Then the statement of \cite{KMOU} reads as follows.

\begin{theorem}[{\rm \cite{KMOU}}] \label{th:iso}
Under the assumptions (I-IV), we have 
\[
\B_{w^{(k)}}(\la)\simeq u_{\la_j}\ot B^{(j)}_a\ot B^{\ot(j-1)}.
\]
\end{theorem}

\subsection{Extremal vectors}
We present a proposition which is convenient to determine 
the extremal vectors. For an element $b$ in $B$ or $\B(\la)$,
let $\ft{i}^{\max}b$ stand for $\ft{i}^{\vphi_i(b)}b$.
Assume $B$ is perfect of level $l$. Therefore, for any $\la
\in\Pcll$ we have the isomorphism.
\begin{equation} \label{eq:iso}
\B(\la)\simeq\P(\la,B).
\end{equation}
For $w\in W$, let us denote by $u_{w\la}$ the vector in the 
lower global 
crystal base of weight $w\la$. Then we have
\begin{eqnarray}
u_{w\la}&=&u_\la\quad\mbox{if }w=1,\nonumber\\
u_{r_iw\la}&=&f_i^{(m)}u_{w\la}\quad
\mbox{if }m=\langle w\la,h_i\rangle\ge0. \label{eq:extre}
\end{eqnarray}
(see section 3.2 of \cite{K2} for this along with the definition of 
$f_i^{(m)}$). note that we can regard $u_{w\la}$ 
as an element of $\b(\la)$ since the $w\la$ weight space is one dimensional. 
again since $r_iw\la$ weight space is one dimensional, in $\b(\la)$, 
using \quad(\ref{eq:extre}) we can and do identify
$u_{r_iw\la}$ with $\ft{i}^m u_{w\la}$. these $u_{w\la}$ are called 
{\em extremal vectors}. We want to find out 
the extremal vector $u_{w^{(k)}\la}$ in the right hand side of 
(\ref{eq:iso}). For any $j\ge1$, we set
\[
b_0^{(j)}=\bbar_j,\qquad 
b_a^{(j)}=\ft{i_a^{(j)}}^{\max}b_{a-1}^{(j)}\quad(a=1,\cd,d).
\]

\begin{proposition} \label{prop:IV'}
Assume that condition (III) holds.
If $\veps_{i_{a+1}^{(j)}}(b_a^{(j)})=0,\linebreak
\vphi_{i_{a+1}^{(j)}}(b_a^{(j)})>0$ $(a=1,\cd,d)$, 
$b_0^{(j+1)}=\ft{i_1^{(j+1)}}^{m^{(j+1)}}b_d^{(j)}$
with $m^{(j+1)}=\langle\la_{j+1},h_{i^{(j+1)}_1}\rangle$ 
for any $j\ge1$, then we have
\[
u_{w^{(k)}\la}=\cd\ot\bbar_{j+2}\ot\bbar_{j+1}\ot(b_a^{(j)})^{\ot j}.
\]
Here, if $k=0$, we set $j=1,a=0$, and otherwise, $j$ and $a$ are 
fixed from $k$ by $k=(j-1)d+a,j\ge1,1\le a\le d$. Also
$i_{d+1}^{(j)}$ should be understood as $i_1^{(j+1)}$.
\end{proposition}

\Proof
We prove by induction on $k$. If $k=0$, the statement should be
understood as
\[
u_{w^{(0)}\la}=\cd\ot\bbar_{3}\ot\bbar_{2}\ot\bbar_{1}.
\]
Thus it is trivial for $w^{(0)}=1$.

Now assume $k>0$. First, we assume $a\ne1$. From the induction 
hypothesis, we have 
\begin{eqnarray*}
u_{w^{(k-1)}\la}&=&\cd\ot\bbar_{j+2}\ot\bbar_{j+1}\ot(b_{a-1}^{(j)})^{\ot j}\\
&\simeq&u_{\la_j}\ot(b_{a-1}^{(j)})^{\ot j}.
\end{eqnarray*}
{}From the assumption of the proposition and the condition (III), we have
\[
\veps_{i^{(j)}_a}(b^{(j)}_{a-1})=0,
\vphi_{i^{(j)}_a}(b^{(j)}_{a-1})>0,
\langle\la_j,h_{i^{(j)}_a}\rangle=0.
\]
Setting $\nu=\vphi_{i^{(j)}_a}(b^{(j)}_{a-1})$, the $i^{(j)}_a$-signature of 
$u_{w^{(k-1)}\la}$ reads as $(\underbrace{+^\nu,\cd,+^\nu}_j)$.
Since $\langle w^{(k-1)}\la,h_{i^{(j)}_a}\rangle
=j(\vphi_{i^{(j)}_a}(b^{(j)}_{a-1})-\veps_{i^{(j)}_a}(b^{(j)}_{a-1}))
=j\nu$, it is now clear that
\begin{eqnarray*}
u_{w^{(k)}\la}&=&u_{r_{i^{(j)}_a}w^{(k-1)}\la}
=\ft{i^{(j)}_a}^{j\nu}u_{w^{(k-1)}\la}\\
&=&\cd\ot\bbar_{j+2}\ot\bbar_{j+1}
\ot(\ft{i^{(j)}_a}^\nu b_{a-1}^{(j)})^{\ot j}\\
&=&\cd\ot\bbar_{j+2}\ot\bbar_{j+1}\ot(b_a^{(j)})^{\ot j}.
\end{eqnarray*}

Taking the assumption $b_0^{(j+1)}=\ft{i_1^{(j+1)}}^{m^{(j+1)}}b_d^{(j)}$
into account, the proof in the case of $a=1$ can be done similarly. 
This completes the proof.
\qed
\medskip

This proposition can be used to check condition (IV). In fact,
assuming the above proposition we can show (IV).
Let us consider the situation where condition (IV) is
satisfied up to $w^{(k)}$. Then
we easily have $\veps_{i_{a+1}^{(j)}}(u_{w^{(k)}\la})=0$,
$\langle w^{(k)}\la,h_{i_{a+1}^{(j)}}\rangle
=\vphi_{i_{a+1}^{(j)}}(u_{w^{(k)}\la})>0$.
Noting Proposition 2 in \cite{KMOU}, we get 
$w^{(k+1)}=r_{i_{a+1}^{(j)}}w^{(k)}\succ w^{(k)}$,
which proves (IV) for $k+1$.
In the cases we shall deal with in this paper, we choose suitable 
known perfect crystals hence condition (I)
already holds. Thus, we are going to check (II), (III), and condition
(IV'):
\begin{description}
\item[(IV')  ]
For any $j\ge1$, $\veps_{i_{a+1}^{(j)}}(b_a^{(j)})=0,
\vphi_{i_{a+1}^{(j)}}(b_a^{(j)})>0$ $(a=1,\cd,d)$ and
$b_0^{(j+1)}=\ft{i_1^{(j+1)}}^{m^{(j+1)}}b_d^{(j)}$ with 
$m^{(j+1)}=\langle\la_{j+1},h_{i^{(j+1)}_1}\rangle$
($i^{(j)}_{d+1}$ should be understood as $i^{(j+1)}_1$.)
\end{description}
which implies condition (IV).

\section{$A_n^{(1)}$ case}

\begin{figure}[ht] \label{fig:A}
\hspace{1cm}
\begin{picture}(300,60)(0,0)
\put(5,20){\circle{10}}
\put(10,20){\line(1,0){40}}
\put(55,20){\circle{10}}
\put(60,20){\line(1,0){40}}
\multiput(110,20)(6,0){10}{\circle*{1}}
\put(173,20){\line(1,0){40}}
\put(218,20){\circle{10}}
\put(223,20){\line(1,0){40}}
\put(268,20){\circle{10}}
\put(0,3){\makebox(10,10){\bf 1}}
\put(50,3){\makebox(10,10){\bf 1}}
\put(213,3){\makebox(10,10){\bf 1}}
\put(263,3){\makebox(10,10){\bf 1}}
\bezier{500}(268,25)(131,80)(5,25)
\end{picture}
\caption{Dynkin diagram $A^{(1)}_n\; (n \ge 2)$}
\end{figure}
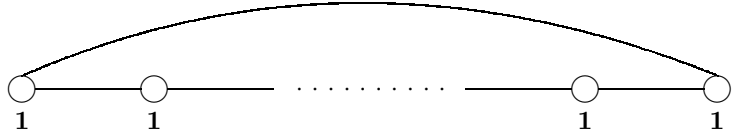

In this section we give explicit descriptions of the $A^{(1)}_{n}$ 
Demazure crystals \linebreak $\B_{w}(l\La_0)$, $l \geq 1$ for a suitably 
chosen linearly ordered chain of Weyl group elements  $w \in 
\{w^{(k)}| k\geq 0\}$. (Note that due to the $\Z_{n+1}$-symmetry we can 
obtain $\B_{w}(l\La_i)$ for all $i$.)
Our starting point is the perfect crystal $B$ 
of level $l$ which is isomorphic to $B(l\ol{\La}_k)$ as crystals 
for $U_{q}(A_{n})$ \cite{KMN2}. As it can be seen the mixing index 
$\kappa =1$ in these cases. The Dynkin diagram $A^{(1)}_n\; (n \ge 2)$
is shown in Figure \ref{fig:A}. The labels are the levels of the 
fundamental weights corresponding to the vertices.

\subsection{Description of the perfect crystal}
Let $I=\Z/(n+1)\Z$ be the index set of the simple roots for
$U_q(A^{(1)}_n)$ and let $J=\{1,\cd,n\}$ be that for $U_q(A_n)$.
For $i\in I$ we define $\iota^{(i)}:J\rightarrow I$ by 
$\iota^{(i)}(j)=i+j$ mod $n+1$, $\ol{\iota}^{(i)}:J\rightarrow I$ by 
$\ol{\iota}^{(i)}(j)=i-j$ mod $n+1$. 
Let $B(l\ol{\La}_k)$ ($1\le k\le n$) be the crystal of type $A_n$
as described in \cite{KN}.
Then for any integers $k,l$ such that $1\le k\le n,l\ge1$, there
exists a unique crystal $B=B^{k,l}$ of type $A^{(1)}_n$ such that
$\iota^{(i)*}(B^{k,l})=B(l\ol{\La}_k)$ and $\ol{\iota}^{(i)*}(B^{k,l})
=B(l\ol{\La}_{k'})$ for all $i$, where $k'=n+1-k$. This $B^{k,l}$ is perfect
of level $l$.

We review the crystal $B(l\ol{\La}_k)$ here. Set $K=\{1,2,\cd,n,n+1\}$. 
With each $b\in B(l\ol{\La}_k)$, we associate a table 
$(m_{i,i'})_{\{1\le i\le k,1\le i'\le l\}}=m(b)$ where 
$m_{i,i'}\in K,m_{i,i'}\le m_{i,i'+1}$ and $m_{i,i'}<m_{i+1,i'}$.
To see the actions of $\et{j}$ and $\ft{j}$ ($j\ne0$), we begin with
the case of $l=1$, each element in $B(\ol{\La}_k)$ corresponds to a column
$(m_i)_{1\le i\le k}$ such that $m_i\in K,m_i<m_{i+1}$. $\et{j}(m_i)$
(resp. $\ft{j}(m_i)$) is obtained by changing $j+1$ in the column to
$j$ (resp. $j$ to $j+1$). If there is no $j+1$ (resp. $j$) in the 
column, or the result breaks the requirements $m_i<m_{i+1}$ ($1\le i<k$),
the action should be set to $0$. To deal with the general case, 
we define the following injection.
\begin{eqnarray}
B(l\ol{\La}_k)&\hookrightarrow&
B(\ol{\La}_k)\ot B(\ol{\La}_k)\ot\cd\ot B(\ol{\La}_k)
\label{eq:inj_A}\\
(m_{i,i'})&\mapsto&(m^{(l)}_i)\ot(m^{(l-1)}_i)\ot\cd\ot(m^{(1)}_i),
\nonumber
\end{eqnarray}
where $(m^{(a)}_i)$ denotes the $a$-th column in the table $(m_{i,i'})$.
Then, the actions on $B(l\ol{\La}_k)$ are defined so that they commute with
this injection. The actions on the right hand side of (\ref{eq:inj_A}) can be
calculated using the signature rule in section \ref{subsec:signature}.

The description of the bijection $\sigma$ is as follows. Let us 
use the notation $\la=(m_0,m_1,\cd,m_n)$ for $\la=\sum_{i=0}^n
m_i\La_i$. Then, 
\[
\sigma:(m_0,m_1,\cd,m_{k'-1},m_{k'},\cd,m_n)\mapsto
(m_k,m_{k+1},\cd,m_n,m_0,\cd,m_{k-1}).
\]
Next we describe $\pbar$.
Fixing $j$, we consider the following sequence of integers:
\[
n+2-jk,n+3-jk,\cd,n+1-(j-1)k.
\]
If there is no integer congruent to $0$ mod $n+1$ in the above sequence,
then $m(\bbar_j)$ is given by $m_{i,i'}\equiv i+n+1-jk$ mod $n+1$. If 
it is not the case, let $\alpha$ be the position of the integer congruent
to $0$ mod $n+1$. Then, $m(\bbar_j)$ is given by $m_{i,i'}=i$ ($1\le i\le
k-\alpha$), $=i+k'$ ($k-\alpha<i\le k$). Here and in what follows in this
section, we set $k'=n+1-k$. In either case, $m_{i,i'}$ does
not depend on $i'$.

\subsection{Description of Demazure crystals}
We set $d=kk'$. We wish to present $i_a^{(j)}$. For this purpose,
it is convenient to associate $(g,r)$ with $a$. Fixing $a$ ($1\le a\le d$),
we set
\[
g=\left[\frac{a-1}{k'}\right],\qquad r=a-1-k'g,
\]
where $[\gamma]$ denotes the largest integer not exceeding $\gamma$.
Note that $0\le g\le k-1,0\le r\le k'-1$. Then, $i_a^{(j)}$ is given by
\[
i^{(j)}_a=k(1-j)-g+r.
\]
Now we have

\begin{theorem}
With the above choice of $i_a^{(j)}$ and $d$, $B$ satisfies (II),(III) and 
(IV').
\end{theorem}

Thanks to the symmetry under $\iota^{(i)*}$, we can reduce checking the
conditions to that for a particular $j$. We choose it to be $n+1$.
Then we have $i_a^{(n+1)}=k-g+r$. Considering the ranges of $g$ and $r$,
we know $i_a^{(n+1)}\not\equiv0$ for any $a$ ($1\le a\le d$). This reduces
the problem for $B$ to that for $B(l\ol{\La}_k)$. We need the following
data: $m(\bbar_{n+1})$ is given by $m_{i,i'}=i$, and $m(\bbar_{n+2})$ 
by $m_{i,i'}=i+k'$.

\begin{proposition}
For $a$, take $g,r$ as above. Then,
\begin{description}
\item[(1)]
$B^{(n+1)}_a=\{m_{i,i'}\mid m_{i,i'}=i\,(i<k-g),m_{k-g,i'}\le k-g+r+1\}$.
\item[(2)]
$m(b_a^{(n+1)})$ is given by $m_{i,i'}=i$ $(i<k-g)$, $=k-g+r+1$ $(i=k-g)$,
$=i+k'$ $(i>k-g)$.
\end{description}
\end{proposition}

\Proof 
We prove by induction on $a$. First, look at (1). If $a=1$, we have
$g=r=0,i^{(n+1)}_1=k$. The assertion is clear from the application
rule through (\ref{eq:inj_A}). Now set $\ol{B}_a$ to be the 
r.h.s. of (1). For any $a$ ($1<a\le d$), we are to show
\[
\ol{B}_a=\bigcup_{m\ge0}\ft{i^{(n+1)}_a}^m
\ol{B}_{a-1}\setminus\{0\},
\]
which is equivalent to show both
\[
\ol{B}_a\supset\bigcup_{m\ge0}\ft{i^{(n+1)}_a}^m
\ol{B}_{a-1}\setminus\{0\} \quad\mbox{and}\quad
\ol{B}_{a-1}\supset\bigcup_{m\ge0}\et{i^{(n+1)}_a}^m
\ol{B}_a\setminus\{0\}.
\]
These are also clear from the application rule.

The proof of (2) is similar.
\qed

\section{$B^{(1)}_{n}$ case} \label{sec:B}

\begin{figure}[ht] \label{fig:B}
\vskip0.7cm
\hspace{0.3cm}
\begin{picture}(300,60)(0,0)
\put(5,20){\circle{10}}
\put(10,20){\line(1,0){40}}
\put(55,20){\circle{10}}
\put(55,25){\line(0,1){40}}
\put(55,70){\circle{10}}
\put(60,20){\line(1,0){40}}
\put(105,20){\circle{10}}
\put(110,20){\line(1,0){40}}
\multiput(160,20)(6,0){10}{\circle*{1}}
\put(223,20){\line(1,0){40}}
\put(268,20){\circle{10}}
\put(273,18){\line(1,0){36}}
\put(273,22){\line(1,0){36}}
\put(313,20){\line(-2,1){10}}
\put(313,20){\line(-2,-1){10}}
\put(318,20){\circle{10}}
\put(62,65){\makebox(10,10){\bf 1}}
\put(0,3){\makebox(10,10){\bf 1}}
\put(50,3){\makebox(10,10){\bf 2}}
\put(100,3){\makebox(10,10){\bf 2}}
\put(263,3){\makebox(10,10){\bf 2}}
\put(313,3){\makebox(10,10){\bf 1}}
\end{picture}
\caption{Dynkin diagram $B^{(1)}_n\; (n \ge 3)$}
\end{figure}
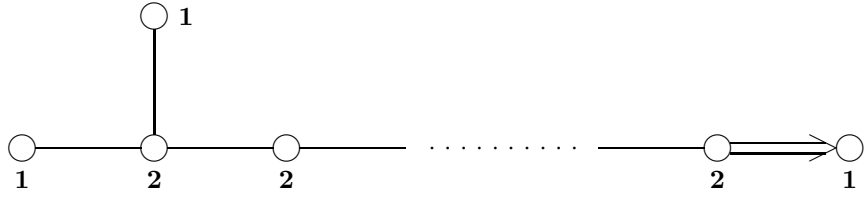

In this section we give explicit descriptions of the $B^{(1)}_{n}$ 
Demazure crystals $\B_{w}(l\La)$, $l \geq 1$, $\La = 
\La_{0},\La_{1},\La_{n}$ (level one dominant weights) for a suitably 
chosen linearly ordered chain of Weyl group elements  $w \in 
\{w^{(k)}| k\geq 0\}$. Our starting point is the perfect crystal $B$ 
of level $l$ which is isomorphic to $B(l\ol{\La}_{1})$ as crystals 
for $U_{q}(B_{n})$ \cite{KMN2}. As it can be seen the mixing index 
$\kappa =1$ in these cases. The Dynkin diagram $B^{(1)}_n\; (n \ge 3)$
is shown in Figure \ref{fig:B}. The labels are the levels of the 
fundamental weights corresponding to the vertices.

\subsection{Description of the perfect crystal}
For any integer $l\geq 1$, we recall the perfect crystal $B$ from 
\cite{KMN2}. As a set
\begin{equation} \label{eq:BofB}
B = \left\{ (x_{1},\cd ,x_{n},x_{0},\ol{x}_{n},\cd ,\ol{x}_{1}) \in 
\Z^{2n} \times\{0,1\} 
\left| 
\begin{array}{c} x_{0} = 0\; \mbox{or} \; 1, x_{i}, \ol{x}_{i} \geq 0, \\ 
\sum_{i=1}^{n}(x_{i}+\ol{x}_{i}) =l 
\end{array}\right.
\right\}.
\end{equation}
The action of $\ft{i}$ is defined as follows: 
For $b=(x_{1},\cd ,x_{n},x_{0},\ol{x}_{n},\cd ,
\ol{x}_{1}) \in B$
\begin{eqnarray*}
\ft{0}b &=& \left\{
\begin{array}{ll}
 (x_{1}, x_{2}+1,\cd ,\ol{x}_{2},\ol{x}_{1}-1) & \mbox{if} \; x_{2}\geq 
 \ol{x}_{2}, \\
 (x_{1}+1, x_{2},\cd ,\ol{x}_{2}-1,\ol{x}_{1}) & \mbox{if} \; x_{2}< 
 \ol{x}_{2},
\end{array}\right. \\
\ft{i}b &=& \left\{
\begin{array}{ll}
 (x_{1},\cd , x_{i}-1,x_{i+1}+1,\cd ,\ol{x}_{1}) & \mbox{if} \; 
 x_{i+1}\geq \ol{x}_{i+1} , \\
 (x_{1},\cd ,\ol{x}_{i+1}-1,\ol{x}_{i}+1,\cd ,\ol{x}_{1}) & \mbox{if} \; 
 x_{i+1} < \ol{x}_{i+1} ,
\end{array}\right. \\
&&\quad \mbox{for $i=1,2,\cdots,n-1$, and}\nonumber \\
\ft{n}b &=& \left\{
\begin{array}{ll}
 (x_{1},\cd , x_{n}-1,x_{0}+1,\ol{x}_{n},\cd ,\ol{x}_{1}) & \mbox{if} \; 
 x_{0}=0 , \\
 (x_{1},\cd , x_{n},x_{0}-1,\ol{x}_{n}+1,\cd ,\ol{x}_{1}) & \mbox{if} \; 
 x_{0} =1 .
\end{array}\right. \\
\end{eqnarray*}
The action of $\et{i}$ is given by $\et{i}b = b'$ if and only if 
$\ft{i}b' = b$. $\vphi_{i}$ and $\veps_{i}$ are given by
\begin{eqnarray*}
&& \vphi_{0}(b) = \ol{x}_{1}+(\ol{x}_{2}-x_{2})_{+},
\quad \vphi_{n}(b) = 2x_{n}+x_{0},\\
&& \vphi_{i}(b) = x_{i}+(\ol{x}_{i+1}-x_{i+1})_{+},
\quad \veps_{i}(b) = \ol{x}_{i}+(x_{i+1}-\ol{x}_{i+1})_{+},\\
&&\quad \mbox{for $i=1,2,\cdots,n-1$, and}\\
&& \veps_{0}(b) = x_{1}+(x_{2}-\ol{x}_{2})_{+},
\quad \veps_{n}(b) = 2\ol{x}_{n}+x_{0},
\end{eqnarray*}
where $(x)_{+}= \max (x,0)$.
Also $\wt(b) = \sum^{n}_{i=0} (\vphi_{i}(b) - \veps_{i}(b))\La_{i}$.

Furthermore, in this case the automorphism $\sigma$ is given by 
\[
\sigma : (m_{0},m_{1},m_{2},\cd ,m_{n}) \longrightarrow 
(m_{1},m_{0},m_{2},\cd ,m_{n})
\]
for $\la = \sum^{n}_{i=0} m_{i}\La_{i}$. Choose $b(l\La_{i}) \in B$ 
as follows:
\[
b(i)=b(l\La_{i}) = \left\{
\begin{array}{ll}
(0,0,\cd,l) & i=0, \\
(l,0,\cd,0) & i=1, \\
(0,\cd ,0,m,x_{0},m,0,\cd ,0) & i=n,
\end{array}\right.
\]
where $x_{0}=\eps(l)$ and $m=\frac{1}{2}(l-x_{0})$.
Here and in what follows, we use the function
\begin{equation} \label{eq:def_eps}
\eps(i)=\left\{
\begin{array}{ll}
0&i:\mbox{even},\\
1&i:\mbox{odd}.
\end{array}
\right.
\end{equation}
Then note that $\vphi(b(i))=l\La_{i},\; i=0,1,n$. So the 
highest weight vector $u_{\la}$ with $\la = l\La_{i}$ corresponds to the 
ground state path $\pbar = \cd\ot\bbar_{k}\ot\cd\ot\bbar_{2}\ot\bbar_{1}$ 
as
\[
\bbar_{k} = \left\{
\begin{array}{ll}
b(\eps(k+i+1)) & i=0,1,\\
b(n) & i=n.
\end{array}\right.
\]
Thus we have, for $\la=l\La_{i},\; i=0,1,n$,
\[
\la_{k} = \left\{
\begin{array}{ll}
l\La_{\eps(k+i)} & i=0 \;or\; 1, \\
l\La_{n} & i=n.
\end{array}\right.
\]

\subsection{Description of Demazure crystals}
For $\la = l\La_{0},\; l\geq 1$, set $d=2n-1$ and
choose the sequence $\{i^{(j)}_{a} | 
j\geq 1, 1\leq a \leq 2n-1\}$ defined as follows:
\begin{eqnarray*}
i^{(j)}_{1} &=& i^{(j)}_{2n-1} = \eps(j+1), \\
i^{(j)}_{a} &=& i^{(j)}_{2n-a} = a \;\; \mbox{for } 2\leq a \leq n. 
\end{eqnarray*}

\begin{theorem}
For $\la = l\La_{0}$, with the above choice of $i^{(j)}_{a}$ and $d$, $B$ 
statisfies (II),(III) and (IV'). Furthermore, 
in this case,
\[
B^{(j)}_{0} = \left\{
\begin{array}{ll}
\{ (0,0,\cd ,l)\} & j\; \mbox{odd}, \\
\{ (l,0,\cd ,0)\} & j\; \mbox{even}, 
\end{array}\right.
\quad B^{(j)}_{2n-1} = B,
\]
and for $1\leq a\leq n-1$, $B^{(j)}_{a}$ and $B^{(j)}_{n+a-1} 
\subseteq B$ are given as follows:
\begin{eqnarray*}
B^{(j)}_{a} &=& \left\{
\begin{array}{ll}
\{ (0,x_{2},\cd ,x_{a+1},0,\cd ,0,\ol{x}_{1})\} & j\; \mbox{odd}, \\
\{ (x_{1},x_{2},\cd ,x_{a+1},0,\cd ,0)\} & j\; \mbox{even}, 
\end{array}\right. \\
B^{(j)}_{n+a-1} &=& \left\{
\begin{array}{ll}
\{ (0,x_{2},\cd ,\ol{x}_{n-a+1},0,\cd ,0,\ol{x}_{1})\} & j\; \mbox{odd}, \\
\{ (x_{1},x_{2},\cd ,\ol{x}_{n-a+1},0,\cd ,0)\} & j\; \mbox{even}.
\end{array}\right. 
\end{eqnarray*}
Here on each set $x_i$'s and $\ol{x}_i$'s run over non negative integers
satisfying the conditions in (\ref{eq:BofB}).
Also, $b^{(j)}_a$ are given as follows:
\[
b^{(j)}_{0} = \left\{
\begin{array}{ll}
\{ (0,0,\cd ,l)\} & j\; \mbox{odd}, \\
\{ (l,0,\cd ,0)\} & j\; \mbox{even}, 
\end{array}\right.
\quad b^{(j)}_{2n-1} = \left\{
\begin{array}{ll}
\{ (l,0,\cd ,0)\} & j\; \mbox{odd}, \\
\{ (0,0,\cd ,l)\} & j\; \mbox{even}, 
\end{array}\right.
\]
and for $1\leq a\leq n-1$, 
\[
b^{(j)}_{a} =
(0,\cd ,0,\displaystyle{\mathop{l}_{\scriptstyle a+1}},0,\cd ,0),
\quad b^{(j)}_{n+a-1} =
(0,\cd ,0,\displaystyle{\mathop{l}_{\scriptstyle \ol{n-a+1}}},0,\cd ,0).
\]
\end{theorem}

\Proof
Since the $j$:odd case is similar, we only consider the $j$:even case, 
where $\la_j=l\La_0$, $\bbar_j=(l,0,\cd,0)$, $i^{(j)}_a=i^{(j)}_{2n-a}=a$
($1\le a\le n$).

Recalling the definitions in section \ref{subsec:crit}, the determination
of the subset $B^{(j)}_a$ for $0\le a\le n$ is easy. Let us consider
$B^{(j)}_{n+1}$. By definition,
\[
B^{(j)}_{n+1}=\bigcup_{m\ge0}\ft{n-1}^m
\{(x_1,\cd,x_n,x_0,\ol{x}_n,0,\cd,0)\}\setminus\{0\}.
\]
In view of the rule of $\ft{n-1}$, we see
\begin{eqnarray*}
&&(x_1,\cd,x_n,x_0,\ol{x}_n,\ol{x}_{n-1},0,\cd,0)\\
&&\quad=\ft{n-1}^{\ol{x}_{n-1}+z_n}
(x_1,\cd,x_{n-1}+z_n,x_n-z_n,x_0,\ol{x}_n+\ol{x}_{n-1},0,\cd,0),
\end{eqnarray*}
where $z_n=(x_n-\ol{x}_n)_+$. This formula proves
\[
B^{(j)}_{n+1}=\{(x_1,\cd,x_n,x_0,\ol{x}_n,\ol{x}_{n-1},0,\cd,0)\}.
\]
The other cases are similar. Thus we have $B^{(j)}_{2n-1}=B$, and
checked the condition (II) with $d=2n-1$. Since
$\langle\la_j,h_{i^{(j)}_a}\rangle=\langle l\La_0,h_{i^{(j)}_a}\rangle=0$
for $1\le a\le d$, the condition (III) is trivial.

The calculation of $b^{(j)}_a$ is simpler. From the rules of
$\veps_i$ and $\vphi_i$, we see $\veps_{i^{(j)}_{a+1}}(b^{(j)}_a)=0$,
$\vphi_{i^{(j)}_{a+1}}(b^{(j)}_a)=l$ ($a\ne n-1$), $2l$ ($a=n-1$), 
$\langle\la_{j+1},h_{i^{(j+1)}_1}\rangle=0$ and $b^{(j+1)}_0=b^{(j)}_d$.
Therefore, we have checked (IV').
\qed

\begin{remark}
In an analogous manner it can be seen, Theorem \ref{th:iso} holds for 
$\la = l\La_{i},\; i=1,n$ if we choose $i^{(i)}_{a}$ as follows:
\begin{eqnarray*}
i = 1 &:& i^{(j)}_{1} = i^{(j)}_{2n-1} = \eps(j), \\
 &&i^{(j)}_{a} = i^{(j)}_{2n-a} = a \;\; for \;2\leq a\leq n. \nonumber\\
i = n &:& i^{(j)}_{a} = n-a+1 \;\;for \;1\leq a\leq n-1, \\
 &&(i^{(j)}_{n}, i^{(j)}_{n+1}) = (1,0)\; or \; (0,1), \nonumber\\
 && i^{(j)}_{n+a} = a \;\;for \;2\leq a\leq n-1. \nonumber 
\end{eqnarray*}
\end{remark}

\section{$D^{(1)}_{n}$ case}

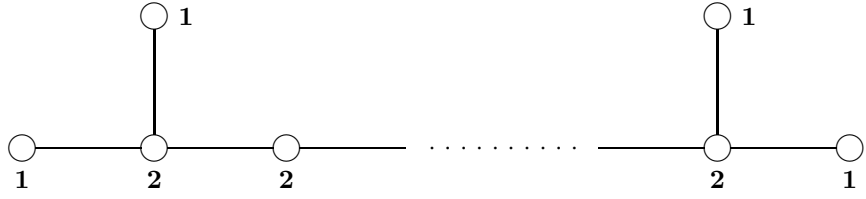
\begin{figure}[ht] \label{fig:D}
\vskip0.7cm
\hspace{0.3cm}
\begin{picture}(300,60)(0,0)
\put(5,20){\circle{10}}
\put(10,20){\line(1,0){40}}
\put(55,20){\circle{10}}
\put(55,25){\line(0,1){40}}
\put(55,70){\circle{10}}
\put(60,20){\line(1,0){40}}
\put(105,20){\circle{10}}
\put(110,20){\line(1,0){40}}
\multiput(160,20)(6,0){10}{\circle*{1}}
\put(223,20){\line(1,0){40}}
\put(268,20){\circle{10}}
\put(268,25){\line(0,1){40}}
\put(268,70){\circle{10}}
\put(273,20){\line(1,0){40}}
\put(318,20){\circle{10}}
\put(62,65){\makebox(10,10){\bf 1}}
\put(0,3){\makebox(10,10){\bf 1}}
\put(50,3){\makebox(10,10){\bf 2}}
\put(100,3){\makebox(10,10){\bf 2}}
\put(263,3){\makebox(10,10){\bf 2}}
\put(275,65){\makebox(10,10){\bf 1}}
\put(313,3){\makebox(10,10){\bf 1}}
\end{picture}
\caption{Dynkin diagram $D^{(1)}_n\; (n \ge 4)$}
\end{figure}

In this section we give explicit descriptions of the $D^{(1)}_{n}$ 
Demazure crystals $\B_{w}(l\La)$, $l \geq 1$, $\La = 
\La_{0},\La_{1},\La_{n-1},\La_{n}$ (level one dominant weights) for a suitably 
chosen linearly ordered chain of Weyl group elements  $w \in 
\{w^{(k)}| k\geq 0\}$. Our starting point is the perfect crystal $B$ 
of level $l$ which is isomorphic to $B(l\ol{\La}_{1})$ as crystals 
for $U_{q}(D_{n})$ \cite{KMN2}. As it can be seen the mixing index 
$\kappa =1$ in these cases. The Dynkin diagram $D^{(1)}_n\; (n \ge 4)$
is shown in Figure \ref{fig:D}. The labels are the levels of the 
fundamental weights corresponding to the vertices.

\subsection{Description of the perfect crystal}
For any integer $l\geq 1$, we recall the perfect crystal $B$ from 
\cite{KMN2}. As a set
\begin{equation} \label{eq:BofD}
B = \left\{ (x_{1},\cd ,x_{n},\ol{x}_{n},\cd ,\ol{x}_{1}) \in 
\Z^{2n}  \left| 
\begin{array}{c}
x_{n} = 0\; \mbox{or} \;\ol{x}_{n} =0,\; x_{i}, \ol{x}_{i} \geq 0,\\ 
\sum_{i=1}^{n}(x_{i}+\ol{x}_{i}) =l 
\end{array}\right.
\right\}. 
\end{equation}
The action of $\ft{i}$ is defined as follows: 
For $b=(x_{1},\cd ,x_{n},\ol{x}_{n},\cd ,\ol{x}_{1}) \in B$
\begin{eqnarray*}
\ft{0}b &=& \left\{
\begin{array}{ll}
 (x_{1}, x_{2}+1,\cd ,\ol{x}_{2},\ol{x}_{1}-1) & \mbox{if} \; x_{2}\geq 
 \ol{x}_{2}, \\
 (x_{1}+1, x_{2},\cd ,\ol{x}_{2}-1,\ol{x}_{1}) & \mbox{if} \; x_{2}< 
 \ol{x}_{2},
\end{array}\right. \\
\ft{i}b &=& \left\{
\begin{array}{ll}
 (x_{1},\cd , x_{i}-1,x_{i+1}+1,\cd ,\ol{x}_{1}) & \mbox{if} \; 
 x_{i+1}\geq \ol{x}_{i+1} , \\
 (x_{1},\cd ,\ol{x}_{i+1}-1,\ol{x}_{i}+1,\cd ,\ol{x}_{1}) & \mbox{if} \; 
 x_{i+1} < \ol{x}_{i+1} ,
\end{array}\right. \\
&&\quad \mbox{for $i=1,2,\cdots,n-2$, and}\\
\ft{n-1}b &=& \left\{
\begin{array}{ll}
 (x_{1},\cd , x_{n-1}-1,x_{n}+1,\cd ,\ol{x}_{1}) & \mbox{if} \; 
 x_{n} \geq 0,\ol{x}_{n}=0 , \\
 (x_{1},\cd ,x_{n},\ol{x}_{n}-1,\ol{x}_{n-1}+1\cd ,\ol{x}_{1}) & \mbox{if} \; 
 x_{n} =0 ,\ol{x}_{n} \geq 1,
\end{array}\right. \\
\ft{n}b &=& \left\{
\begin{array}{ll}
 (x_{1},\cd ,x_{n}-1,\ol{x}_{n},\ol{x}_{n-1}+1,\cd ,\ol{x}_{1}) & \mbox{if} \; 
 x_{n} \geq 1,\ol{x}_{n}=0 , \\
 (x_{1},\cd ,x_{n-1}-1,x_{n},\ol{x}_{n}+1\cd ,\ol{x}_{1}) & \mbox{if} \; 
 x_{n} =0 ,\ol{x}_{n} \geq 0.
\end{array}\right. 
\end{eqnarray*}
The action of $\et{i}$ is given by $\et{i}b = b'$ if and only if 
$\ft{i}b' = b$. $\vphi_{i}$ and $\veps_{i}$ are given by
\begin{eqnarray*}
&& \vphi_{0}(b) = \ol{x}_{1}+(\ol{x}_{2}-x_{2})_{+}, \quad
\vphi_{n-1}(b) = x_{n-1}+\ol{x}_{n}, \quad
\vphi_{n}(b) = x_{n-1}+x_{n},\\
&&\vphi_{i}(b) = x_{i}+(\ol{x}_{i+1}-x_{i+1})_{+}, \quad
\veps_{i}(b) = \ol{x}_{i}+(x_{i+1}-\ol{x}_{i+1})_{+}\\
&&\quad \mbox{for $i=1,2,\cd,n-2$, and}\\ 
&&\veps_{0}(b) = x_{1}+(x_{2}-\ol{x}_{2})_{+}, \quad
\veps_{n-1}(b) = \ol{x}_{n-1}+x_{n}, \quad
\veps_{n}(b) = \ol{x}_{n-1}+\ol{x}_{n},
\end{eqnarray*}
where $(x)_{+}= \max (x,0)$.
Also $\wt(b) = \sum^{n}_{i=0}(\vphi_{i}(b) - \veps_{i}(b))\La_{i}$.

Furthermore, the automorphism $\sigma$ is given by 
\[
\sigma : (m_{0},m_{1},m_{2},\cd ,m_{n-2},m_{n-1},m_{n}) \longrightarrow 
(m_{1},m_{0},m_{2},\cd ,m_{n-2},m_{n},m_{n-1})
\]
for $\la = \sum^{n}_{i=0} m_{i}\La_{i}$. Choose $b(l\La_{i}) \in B$, 
$i=0,1,n-1,n$ as follows:
\[
b(i)=b(l\La_{i}) = \left\{
\begin{array}{ll}
(0,0,\cd,l) & i=0, \\
(l,0,\cd,0) & i=1, \\
(0,\cd ,0,\displaystyle{\mathop{l}_{\scriptstyle \ol{n}}},0,\cd ,0) & i=n-1,\\
(0,\cd ,0,\displaystyle{\mathop{l}_{\scriptstyle n}},0,\cd ,0) & i=n.
\end{array}\right.
\]
Note that $\vphi(b(i))=l\La_{i},\; i=0,1,n-1,n$. the 
highest weight vector $u_{\la}$ with $\la = l\La_{i}$ corresponds to the 
ground state path $\pbar = \cd\ot\bbar_{k}\ot\cd\ot\bbar_{1}$ 
with
\[
\bbar_{k} = \left\{
\begin{array}{ll}
b(\eps(k+i+1)) & i=0,1,\\
b(n-\eps(k+n-i+1)) & i=n-1,n.
\end{array}\right.
\]
Here $\eps(i)$ is defined in (\ref{eq:def_eps}).
Thus we have, for $\la=l\La_{i},\; i=0,1,n-1,n$,
\[
\la_{k} = \left\{
\begin{array}{ll}
l\La_{\eps(k+i)} & i=0,1, \\
l\La_{n-\eps(k+n-i)} & i=n-1,n.
\end{array}\right.
\]

\subsection{Description of Demazure crystals}
For $\la = l\La_{0},\; l\geq 1$, set $d=2n-2$ and
choose the sequence $\{i^{(j)}_{a} | 
j\geq 1, 1\leq a \leq 2n-2\}$ defined as follows:
\begin{eqnarray*}
i^{(j)}_{1} &=& i^{(j)}_{2n-2} = \eps(j+1), \\
(i^{(j)}_{n-1},i^{(j)}_{n}) &=& (n-1,n),\\
i^{(j)}_{a} &=& i^{(j)}_{2n-a-1} = a \;\; for\; 2\leq a \leq n-2. 
\end{eqnarray*}
\begin{theorem}
For $\la = l\La_{0}$, with the above choice of $i^{(j)}_{a}$ and $d$,
$B$ statisfies (II),(III) and (IV'). 
Furthermore, in this case,
\[
B^{(j)}_{0} = \left\{
\begin{array}{ll}
\{ (0,0,\cd ,l)\} & j\; \mbox{odd}, \\
\{ (l,0,\cd ,0)\} & j\; \mbox{even}, 
\end{array}\right.
\quad , B^{(j)}_{2n-2} = B,
\]
and for $1\leq a\leq 2n-3$, $B^{(j)}_{a} \subseteq B$ are given as follows:
\begin{eqnarray*}
B^{(j)}_{a} &=& \left\{
\begin{array}{ll}
\{ (0,x_{2},\cd ,x_{a+1},0,\cd ,0,\ol{x}_{1})\} & j\; \mbox{odd}, \\
\{ (x_{1},x_{2},\cd ,x_{a+1},0,\cd ,0)\} & j\; \mbox{even}, 
\end{array}\right. \\
B^{(j)}_{n-1} &=& \left\{
\begin{array}{ll}
\{ (0,x_{2},\cd ,x_{n},0,\cd ,0,\ol{x}_{1})\} & j\; \mbox{odd}, \\
\{ (x_{1},x_{2},\cd ,x_{n},0,\cd ,0)\} & j\; \mbox{even}, 
\end{array}\right. \\
B^{(j)}_{n+a-1} &=& \left\{
\begin{array}{ll}
\{ (0,x_{2},\cd ,\ol{x}_{n-a},0,\cd ,0,\ol{x}_{1})\} & j\; \mbox{odd}, \\
\{ (x_{1},x_{2},\cd ,\ol{x}_{n-a},0,\cd ,0)\} & j\; \mbox{even}, 
\end{array}\right. 
\end{eqnarray*}
where $1\leq a \leq n-2$.
Here on each set $x_i$'s and $\ol{x}_i$'s run over non negative integers
satisfying the conditions in (\ref{eq:BofD}).
Also, $b^{(j)}_a$ are given as follows:
\[
b^{(j)}_{0} = \left\{
\begin{array}{ll}
\{ (0,0,\cd ,l)\} & j\; \mbox{odd}, \\
\{ (l,0,\cd ,0)\} & j\; \mbox{even}, 
\end{array}\right.
\quad b^{(j)}_{2n-1} = \left\{
\begin{array}{ll}
\{ (l,0,\cd ,0)\} & j\; \mbox{odd}, \\
\{ (0,0,\cd ,l)\} & j\; \mbox{even}, 
\end{array}\right.
\]
\[
b^{(j)}_{n-1} = 
(0,\cd ,0,\displaystyle{\mathop{l}_{\scriptstyle n}},0,\cd ,0),
\]
and for $1\leq a\leq n-2$, 
\[
b^{(j)}_{a} =
(0,\cd ,0,\displaystyle{\mathop{l}_{\scriptstyle a+1}},0,\cd ,0),
\quad b^{(j)}_{n+a-1} =
(0,\cd ,0,\displaystyle{\mathop{l}_{\scriptstyle \ol{n-a}}},0,\cd ,0).
\]
\end{theorem}

\Proof
The proof is similar to $B^{(1)}_n$.
\qed

\begin{remark}
We note that for $\la = l\La_{0}$ choosing 
$(i^{(j)}_{n-1},i^{(j)}_{n}) = (n,n-1)$ also works. Furthermore, for 
$\la = l\La_{i},\; i=1,n-1,n$ it can be seen, Theorem \ref{th:iso} holds for 
the following choices of the sequence $\{i^{(j)}_{a}\}$:
\begin{eqnarray*}
i = 1 &:& i^{(j)}_{1} = i^{(j)}_{2n-2} = \eps(j), \\
&&(i^{(j)}_{n-1},i^{(j)}_{n}) = (n-1,n)\; or \;(n,n-1), \nonumber\\ 
&& i^{(j)}_{a} = i^{(j)}_{2n-a-1} = a \;\; for \;2\leq a\leq n-2, \nonumber\\
i = n-1 &:& i^{(j)}_{1} = i^{(j)}_{2n-2} = n-\eps(j), \\
&&(i^{(j)}_{n-1}, i^{(j)}_{n}) = (1,0)\; or \; (0,1), \nonumber\\
&&i^{(j)}_{a} = i^{(j)}_{2n-a-1} = n-a \;\;for \;2\leq a\leq n-2,\nonumber\\
i = n &:& i^{(j)}_{1} = i^{(j)}_{2n-2} = n-\eps(j+1), \\
&& (i^{(j)}_{n-1}, i^{(j)}_{n}) = (1,0)\; or \; (0,1), \nonumber\\ 
&& i^{(j)}_{a} = i^{(j)}_{2n-a-1} = n-a \;\;for \;2\leq a\leq n-2.\nonumber
\end{eqnarray*}
\end{remark}

\section{$A^{(2)}_{2n-1}$ case}

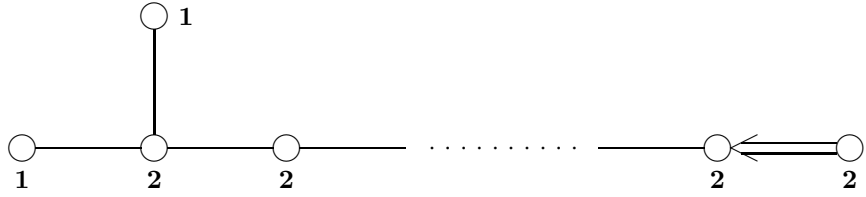
\begin{figure}[ht] \label{fig:A2o}
\vskip0.7cm
\hspace{0.3cm}
\begin{picture}(300,60)(0,0)
\put(5,20){\circle{10}}
\put(10,20){\line(1,0){40}}
\put(55,20){\circle{10}}
\put(55,25){\line(0,1){40}}
\put(55,70){\circle{10}}
\put(60,20){\line(1,0){40}}
\put(105,20){\circle{10}}
\put(110,20){\line(1,0){40}}
\multiput(160,20)(6,0){10}{\circle*{1}}
\put(223,20){\line(1,0){40}}
\put(268,20){\circle{10}}
\put(277,18){\line(1,0){36}}
\put(277,22){\line(1,0){36}}
\put(273,20){\line(2,1){10}}
\put(273,20){\line(2,-1){10}}
\put(318,20){\circle{10}}
\put(62,65){\makebox(10,10){\bf 1}}
\put(0,3){\makebox(10,10){\bf 1}}
\put(50,3){\makebox(10,10){\bf 2}}
\put(100,3){\makebox(10,10){\bf 2}}
\put(263,3){\makebox(10,10){\bf 2}}
\put(313,3){\makebox(10,10){\bf 2}}
\end{picture}
\caption{Dynkin diagram $A^{(2)}_{2n-1}\; (n \ge 3)$}
\end{figure}

In this section we give explicit descriptions of the $A^{(2)}_{2n-1}$ 
Demazure crystals $\B_{w}(l\La)$, $l \geq 1$, $\La = 
\La_{0},\La_{1}$ (level one dominant weights) for a suitably 
chosen linearly ordered chain of Weyl group elements  $w \in 
\{w^{(k)}| k\geq 0\}$. Our starting point is the perfect crystal $B$ 
of level $l$ which is isomorphic to $B(l\ol{\La}_{1})$ as crystals 
for $U_{q}(C_{n})$ \cite{KMN2}. As it can be seen the mixing index 
$\kappa =1$ in these cases. The Dynkin diagram $A^{(2)}_{2n-1}\; (n \ge 3)$
is shown in Figure \ref{fig:A2o}. The labels are the levels of the 
fundamental weights corresponding to the vertices.

\subsection{Description of the perfect crystal}
For any integer $l\geq 1$, we recall the perfect crystal $B$ from 
\cite{KMN2}. As a set
\begin{equation} \label{eq:BofA2o}
B = \{ (x_{1},\cd ,x_{n},\ol{x}_{n},\cd ,\ol{x}_{1}) \in 
\Z^{2n}  | x_{i}, \ol{x}_{i} \geq 0,\sum_{i=1}^{n}(x_{i}+\ol{x}_{i}) =l \}.
\nonumber
\end{equation}
The action of $\ft{i}$ on $B$ are defined as follows: 
For $b=(x_{1},\cd ,x_{n},\ol{x}_{n},\cd ,\ol{x}_{1}) \in B$
\begin{eqnarray*}
\ft{0}b &=& \left\{
\begin{array}{ll}
 (x_{1}, x_{2}+1,\cd ,\ol{x}_{2},\ol{x}_{1}-1) & \mbox{if} \; x_{2}\geq 
 \ol{x}_{2}, \\
 (x_{1}+1, x_{2},\cd ,\ol{x}_{2}-1,\ol{x}_{1}) & \mbox{if} \; x_{2}< 
 \ol{x}_{2},
\end{array}\right. \\
\ft{i}b &=& \left\{
\begin{array}{ll}
 (x_{1},\cd , x_{i}-1,x_{i+1}+1,\cd ,\ol{x}_{1}) & \mbox{if} \; 
 x_{i+1}\geq \ol{x}_{i+1} , \\
 (x_{1},\cd ,\ol{x}_{i+1}-1,\ol{x}_{i}+1,\cd ,\ol{x}_{1}) & \mbox{if} \; 
 x_{i+1} < \ol{x}_{i+1} ,
\end{array}\right. \\
&&\quad \mbox{for $i=1,2,\cd ,n-1$, and}\\
\ft{n}b &=& 
 (x_{1},\cd , x_{n}-1,\ol{x}_{n}+1,\cd ,\ol{x}_{1}). \nonumber
\end{eqnarray*}
The action of $\et{i}$ is given by $\et{i}b = b'$ if and only if 
$\ft{i}b' = b$. $\vphi_{i}$ and $\veps_{i}$ are given by
\begin{eqnarray*}
&&\vphi_{0}(b) = \ol{x}_{1}+(\ol{x}_{2}-x_{2})_{+}, \quad
\vphi_{n}(b) = x_{n},\\
&&\vphi_{i}(b) = x_{i}+(\ol{x}_{i+1}-x_{i+1})_{+}, \quad
\veps_{i}(b) = \ol{x}_{i}+(x_{i+1}-\ol{x}_{i+1})_{+},\\
&&\quad \mbox{for $i=1,2,\cd,n-2$, and}\\ 
&&\veps_{0}(b) = x_{1}+(x_{2}-\ol{x}_{2})_{+}, \quad
\veps_{n}(b) = \ol{x}_{n},
\end{eqnarray*}
where $(x)_+=\max(x,0)$.
Also $\wt(b) = \sum^{n}_{i=0}(\vphi_{i}(b) - \veps_{i}(b))\La_{i}$.

In this case the automorphism $\sigma$ is given by 
\[
\sigma : (m_{0},m_{1},m_{2},\cd ,m_{n}) \longrightarrow 
(m_{1},m_{0},m_{2},\cd ,m_{n})
\]
for $\la = \sum^{n}_{i=0} m_{i}\La_{i}$. Choose $b(i)=b(l\La_{i}) \in B$, 
$i=0,1$ as follows:
\[
b(i)= \left\{
\begin{array}{ll}
(0,0,\cd,l) & i=0, \\
(l,0,\cd,0) & i=1. 
\end{array}\right.
\]
Note that $\vphi(b(i))=l\La_{i},\; i=0,1$. The 
highest weight vector $u_{\la}$ with $\la = l\La_{i}$ corresponds to the 
ground state path $\pbar = \cd\ot\bbar_{k}\ot\cd\ot\bbar_{1}$ 
with
\[
\bbar_{k} = b(\eps(k+i+1)).
\]
Here $\eps(i)$ is defined in (\ref{eq:def_eps}).
We also have, for $\la=l\La_{i},\; i=0,1$,
\[
\la_{k} = l\La_{\eps(k+i)}.
\]

\subsection{Description of Demazure crystals}
For $\la = l\La_{0},\; l\geq 1$, set $d=2n-1$ and
choose the sequence $\{i^{(j)}_{a} | 
j\geq 1, 1\leq a \leq 2n-1\}$ defined as follows:
\begin{eqnarray*}
i^{(j)}_{1} &=& i^{(j)}_{2n-1} = \eps(j+1), \\
i^{(j)}_{a} &=& i^{(j)}_{2n-a} = a \;\; \mbox{for}\; 2\leq a \leq n. 
\end{eqnarray*}
\begin{theorem}
For $\la = l\La_{0}$, with the above choice of $i^{(j)}_{a}$ and $d$, 
$B$ statisfies (II),(III) and (IV'). Furthermore, in this case,
\[
B^{(j)}_{0} = \left\{
\begin{array}{ll}
\{ (0,0,\cd ,l)\} & j\; \mbox{odd}, \\
\{ (l,0,\cd ,0)\} & j\; \mbox{even}, 
\end{array}\right.
\quad B^{(j)}_{2n-1} = B,
\]
and for $1\leq a\leq 2n-2$, $B^{(j)}_{a} \subseteq B$ are given as follows:
\begin{eqnarray*}
B^{(j)}_{a} &=& \left\{
\begin{array}{ll}
\{ (0,x_{2},\cd ,x_{a+1},0,\cd ,0,\ol{x}_{1})\} & j\; \mbox{odd}, \\
\{ (x_{1},x_{2},\cd ,x_{a+1},0,\cd ,0)\} & j\; \mbox{even}, 
\end{array}\right. \\
B^{(j)}_{n+a-1} &=& \left\{
\begin{array}{ll}
\{ (0,x_{2},\cd ,\ol{x}_{n-a+1},0,\cd ,0,\ol{x}_{1})\} & j\; \mbox{odd}, \\
\{ (x_{1},x_{2},\cd ,\ol{x}_{n-a+1},0,\cd ,0)\} & j\; \mbox{even}, 
\end{array}\right. 
\end{eqnarray*}
where $1\leq a \leq n-1$.
Here on each set $x_i$'s and $\ol{x}_i$'s run over non negative integers
satisfying the conditions in (\ref{eq:BofA2o}).
Also, $b^{(j)}_a$ are given as follows:
\[
b^{(j)}_{0} = \left\{
\begin{array}{ll}
\{ (0,0,\cd ,l)\} & j\; \mbox{odd}, \\
\{ (l,0,\cd ,0)\} & j\; \mbox{even}, 
\end{array}\right.
\quad b^{(j)}_{2n-1} = \left\{
\begin{array}{ll}
\{ (l,0,\cd ,0)\} & j\; \mbox{odd}, \\
\{ (0,0,\cd ,l)\} & j\; \mbox{even}, 
\end{array}\right.
\]
and for $1\leq a\leq n-1$, 
\[
b^{(j)}_{a} =
(0,\cd ,0,\displaystyle{\mathop{l}_{\scriptstyle a+1}},0,\cd ,0),
\quad b^{(j)}_{n+a-1} =
(0,\cd ,0,\displaystyle{\mathop{l}_{\scriptstyle \ol{n-a+1}}},0,\cd ,0).
\]
\end{theorem}

\Proof
The proof is similar to $B^{(1)}_n$.
\qed

\begin{remark}
For $\la = l\La_{1}$ it can be seen similarly that Theorem \ref{th:iso}
holds for the following choices of the sequence $\{i^{(j)}_{a}\}$:
\begin{eqnarray*}
i = 1 &:& i^{(j)}_{1} = i^{(j)}_{2n-1} = \eps(j), \\ 
&& i^{(j)}_{a} = i^{(j)}_{2n-a} = a \;\; for \;2\leq a\leq n. \nonumber\\
\end{eqnarray*}
\end{remark}

\section{$A_{2n}^{(2)}$ case}

\begin{figure}[ht] \label{fig:A2e}
\hspace{1cm}
\begin{picture}(300,60)(0,0)
\put(5,20){\circle{10}}
\put(14,18){\line(1,0){36}}
\put(14,22){\line(1,0){36}}
\put(10,20){\line(2,1){10}}
\put(10,20){\line(2,-1){10}}
\put(55,20){\circle{10}}
\put(60,20){\line(1,0){40}}
\multiput(110,20)(6,0){10}{\circle*{1}}
\put(173,20){\line(1,0){40}}
\put(218,20){\circle{10}}
\put(227,18){\line(1,0){36}}
\put(227,22){\line(1,0){36}}
\put(223,20){\line(2,1){10}}
\put(223,20){\line(2,-1){10}}
\put(268,20){\circle{10}}
\put(0,3){\makebox(10,10){\bf 1}}
\put(50,3){\makebox(10,10){\bf 2}}
\put(213,3){\makebox(10,10){\bf 2}}
\put(263,3){\makebox(10,10){\bf 2}}
\end{picture}
\caption{Dynkin diagram $A^{(2)}_{2n}\; (n \ge 2)$}
\end{figure}
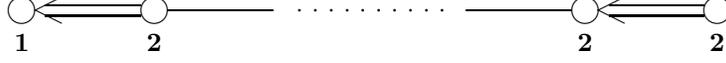

In this section we give explicit descriptions of the $A_{2n}^{(2)}$
Demazure crystals \linebreak $\B_{w}(l\La_0)$, $l \ge 1$ ($\La_0$ being the 
level one dominant weight) for a suitably chosen linearly ordered chain of 
Weyl group elements $w \in \{ w^{(k)} | k \ge 0 \}$.
Our starting point is the perfect crystal $B$ of level $l$ which is
isomorphic to $B(0) \oplus B(\ol{\La}_1) \oplus 
\cdots \oplus B(l \ol{\La}_1)$ as crystals for $U_q(C_n)$ \cite{KMN2}.
It is shown that the mixing index $\kappa = 1$ in these cases.
The Dynkin diagram $A^{(2)}_{2n}\; (n \ge 2)$ is shown in Figure \ref{fig:A2e}.
The labels are the levels of the fundamental weights corresponding to the 
vertices.

\subsection{Description of the perfect crystal}
For any integer $l \ge 1$, we recall the perfect crystal $B$ from \cite{KMN2}.
As a set,
\begin{equation} \label{eq:BofA2e}
B = \{  (x_1,\cdots,x_n,\ol{x}_n,\cdots,\ol{x}_1) \in
\Z^{2n} \;| \; x_i,\ol{x}_i \ge 0, \sum_{i=1}^n (x_i + \ol{x}_i)
\le l \}.
\end{equation}
The action of $\ft{i}$ on $B$ is defined as follows:
For $b=(x_1,\cdots,x_n,\ol{x}_n,
\cdots,\ol{x}_1) \in B$,
\begin{eqnarray*}
&& \ft{0} b = \left\{ 
\begin{array}{ll}
(x_1 +1,x_2,\cdots,\ol{x}_1) & \mbox{if $x_1 \ge \ol{x}_1$},\\
(x_1,\cdots,\ol{x}_2,\ol{x}_1-1) & \mbox{if $x_1 < \ol{x}_1$},
\end{array}
\right.
\nonumber \\
&& \ft{i} b = \left\{ 
\begin{array}{ll}
(x_1,\cdots,x_i-1,x_{i+1}+1,\cdots,\ol{x}_1) & 
\mbox{if $x_{i+1} \ge \ol{x}_{i+1}$},\\
(x_1,\cdots,\ol{x}_{i+1}-1,\ol{x}_i+1,\cdots,\ol{x}_1) & 
\mbox{if $x_{i+1} < \ol{x}_{i+1}$},
\end{array}
\right.
\nonumber \\
&&\quad \mbox{for $i=1,2,\cdots,n-1$, and}\nonumber \\
&& \ft{n} b =(x_1,\cdots,x_n-1,\ol{x}_n+1,\cdots,\ol{x}_1).
\end{eqnarray*}
The action of $\et{i}$ is given by $\et{i}b=b'$ if and only if
$\ft{i}b'=b$. $\vphi_{i}$ and $\veps_{i}$ are given by
\begin{eqnarray*}
&& \vphi_0(b) = l - \sum_{i=1}^n (x_i + \ol{x}_i) + 
2(\ol{x}_1 - x_1)_+,\quad
\vphi_n(b) = x_n ,\nonumber\\
&& \vphi_i(b) = x_i + (\ol{x}_{i+1} - x_{i+1})_+,\quad
\veps_i(b) = \ol{x}_i + (x_{i+1} - \ol{x}_{i+1})_+,\nonumber\\
&&\quad \mbox{for $i=1,2,\cdots,n-1$, and}\nonumber\\
&& \veps_0(b) = l - \sum_{i=1}^n (x_i + \ol{x}_i) + 
2(x_1 - \ol{x}_1)_+,\quad
\veps_n(b) = \ol{x}_n .
\end{eqnarray*}
As before, $(x)_+=\max(x,0)$ and 
$\wt(b)=\sum_{i=0}^n(\vphi_i(b)-\veps_i(b))\La_i$.
In this case the automorphism $\sigma$ is given by
\[
\sigma : (m_0,m_1,\cdots,m_n) \mapsto (m_0,m_1,\cdots,m_n),
\]
for $\la = \sum_{i=0}^n m_i \La_i$.

Choose $b(0) = (0,0,\cdots,0) \in B$.
Then $\vphi(b(0))=l \La_0$.
The highest weight vector $u_\la$ with $\la=l\La_0$ corresponds to the ground
state path $\pbar = \cd \ot \bbar_k \ot \cd \ot \bbar_1$,
where $\bbar_k = b(0)$ for all $k$.
We also have $\la_k = l \La_0$ for all $k$.

\subsection{Description of Demazure crystals}
For $\la = l \La_0,l \ge 1$, set $d=2n$ and choose the sequence $\{i_a^{(j)} |
j \ge 1,1 \le a \le 2n \}$ defined as follows:
\[
i_a^{(j)} = \left\{\begin{array}{ll}
a-1, & 1 \le a \le n+1,\\
2n+1-a, & n+2 \le a \le 2n.
\end{array}
\right.
\]

\begin{theorem}
For $\la = l \La_0$, with the above choice of $i_a^{(j)}$ and $d$,
$B$ satisfies (II),(III) and (IV').
Furthermore, in this case,
\begin{eqnarray*}
&& B_0^{(j)} = \{ (0,\cd ,0) \}, \quad B_{2n}^{(j)} = B, \quad \mbox{and}
\nonumber\\
&& B_a^{(j)} = \{ (x_1,\cd ,x_a,0,\cd,0) \} \subseteq B, \quad 1 \le a
\le n,
\nonumber\\
&& B_{n+a}^{(j)} = \{ (x_1,\cd ,x_n,\ol{x}_n,\cd,\ol{x}_{n-a+1},
0,\cd ,0) \} \subseteq B, \quad 1 \le a
\le n-1.
\nonumber\\
&&
\end{eqnarray*}
Here on each set $x_i$'s and $\ol{x}_i$'s run over non negative integers
satisfying the conditions in (\ref{eq:BofA2e}).
Also, $b^{(j)}_a$ are given as follows:
\[
b^{(j)}_{0} = (0,\cd,0),
\]
and for $1\leq a\leq n$, 
\[
b^{(j)}_{a} =
(0,\cd ,0,\displaystyle{\mathop{l}_{\scriptstyle a}},0,\cd ,0),
\quad b^{(j)}_{n+a} =
(0,\cd ,0,\displaystyle{\mathop{l}_{\scriptstyle \ol{n-a+1}}},0,\cd ,0).
\]
\end{theorem}

\Proof
For any $j$, we have $\la_j=l\La_0$, $\bbar_j=(0,\cd,0)$,
$i^{(j)}_a=a-1$, $i^{(j)}_{n+a}=n-a+1$ ($1\le a\le n$). From the rule
of $\ft{i}$, the determination of the subset $B^{(j)}_a$ for 
$0\le a\le n+1$ is easy. As $B^{(j)}_a$ for $n+2\le a\le 2n$ concerns,
the determination is just the same as $B^{(1)}_n$ in section \ref{sec:B}.
Therefore, the condition (II) is valid with $d=2n$. We have 
$\langle\la_j,h_{i^{(j)}_a}\rangle=0$ ($a\ne1$), $=l$ ($a=1$).
Since $\veps_0((0,\cd,0))=l$, (III) is also valid.

The calculation of $b^{(j)}_a$ is simpler. From the rules of
$\veps_i$ and $\vphi_i$, we see $\veps_{i^{(j)}_{a+1}}(b^{(j)}_a)=0$,
$\vphi_{i^{(j)}_{a+1}}(b^{(j)}_a)=l$ ($a\ne d$), $2l$ ($a=d$), 
$\langle\la_{j+1},h_{i^{(j+1)}_1}\rangle=l$ and 
$b^{(j+1)}_0=\ft{0}^l b^{(j)}_d$.
Therefore, we have checked (IV').
\qed

\section{$D_{n+1}^{(2)}$ case}

\begin{figure}[ht] \label{fig:D2}
\hspace{1cm}
\begin{picture}(300,60)(0,0)
\put(5,20){\circle{10}}
\put(14,18){\line(1,0){36}}
\put(14,22){\line(1,0){36}}
\put(10,20){\line(2,1){10}}
\put(10,20){\line(2,-1){10}}
\put(55,20){\circle{10}}
\put(60,20){\line(1,0){40}}
\multiput(110,20)(6,0){10}{\circle*{1}}
\put(173,20){\line(1,0){40}}
\put(218,20){\circle{10}}
\put(223,18){\line(1,0){36}}
\put(223,22){\line(1,0){36}}
\put(263,20){\line(-2,1){10}}
\put(263,20){\line(-2,-1){10}}
\put(268,20){\circle{10}}
\put(0,3){\makebox(10,10){\bf 1}}
\put(50,3){\makebox(10,10){\bf 2}}
\put(213,3){\makebox(10,10){\bf 2}}
\put(263,3){\makebox(10,10){\bf 1}}
\end{picture}
\caption{Dynkin diagram $D^{(2)}_{n+1}\; (n \ge 2)$}
\end{figure}
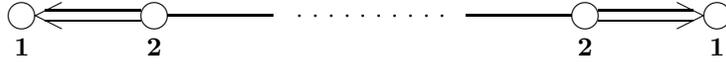

In this section we give explicit descriptions of the $D_{n+1}^{(2)}$
Demazure crystals $\B_{w}(l\La)$, $l \ge 1$, $\La = \La_0,\La_n$
(level one dominant weights) for a suitably chosen linearly ordered 
chain of Weyl group elements $w \in \{ w^{(k)} | k \ge 0 \}$.
Our starting point is the perfect crystal $B$ of level $l$ which is
isomorphic to $B(0) \oplus B(\ol{\La}_1) \oplus 
\cdots \oplus B(l \ol{\La}_1)$ as crystals for $U_q(B_n)$ \cite{KMN2}.
It is shown that the mixing index $\kappa = 1$ in these cases.
The Dynkin diagram $D^{(2)}_{n+1}\; (n \ge 2)$ is shown in Figure \ref{fig:D2}.
The labels are the levels of the fundamental weights corresponding to the 
vertices.

\subsection{Description of the perfect crystal}
For any integer $l \ge 1$, we briefly recall the perfect crystal $B$ from \cite{KMN2}.
As a set,
\begin{equation} \label{eq:BofD2}
B = \left\{  (x_1,\cdots,x_n,x_0,\ol{x}_n,\cdots,\ol{x}_1) \in
\Z^{2n} \times \{0,1\} \Biggm| 
\begin{array}{l}
x_0=\mbox{$0$ or $1$},x_i,\ol{x}_i \ge 0, \\
\sum_{i=1}^n (x_i + \ol{x}_i) \le l 
\end{array}
\right\}.
\end{equation}
The action of $\ft{i}$ on $B$ is defined as follows:
For $b=(x_1,\cdots,x_n,x_0,\ol{x}_n,
\cdots,\ol{x}_1) \in B$,
\begin{eqnarray*}
&& \ft{0} b = \left\{ 
\begin{array}{ll}
(x_1 +1,x_2,\cdots,\ol{x}_1) & \mbox{if $x_1 \ge \ol{x}_1$},\\
(x_1,\cdots,\ol{x}_2,\ol{x}_1-1) & \mbox{if $x_1 < \ol{x}_1$},
\end{array}
\right.
\nonumber \\
&& \ft{i} b = \left\{ 
\begin{array}{ll}
(x_1,\cdots,x_i-1,x_{i+1}+1,\cdots,\ol{x}_1) & 
\mbox{if $x_{i+1} \ge \ol{x}_{i+1}$},\\
(x_1,\cdots,\ol{x}_{i+1}-1,\ol{x}_i+1,\cdots,\ol{x}_1) & 
\mbox{if $x_{i+1} < \ol{x}_{i+1}$},
\end{array}
\right.
\nonumber \\
&&\quad \mbox{for $i=1,2,\cdots,n-1$, and}\nonumber \\
&& \ft{n} b = \left\{ 
\begin{array}{ll}
(x_1,\cdots,x_n-1,x_0+1,\ol{x}_n,\cdots,\ol{x}_1) & 
\mbox{if $x_0 = 0$},\\
(x_1,\cdots,x_n,x_0-1,\ol{x}_n+1,\cdots,\ol{x}_1) & 
\mbox{if $x_0 = 1$}.
\end{array}
\right.
\end{eqnarray*}
The action of $\et{i}$ is given by $\et{i}b=b'$ if and only if
$\ft{i}b'=b$. $\vphi_{i}$ and $\veps_{i}$ are given by
\begin{eqnarray*}
&& \vphi_0(b) = l - x_0 - \sum_{i=1}^n (x_i + \ol{x}_i) + 
2(\ol{x}_1 - x_1)_+,\quad
\vphi_n(b) = 2 x_n + x_0,\nonumber\\
&& \vphi_i(b) = x_i + (\ol{x}_{i+1} - x_{i+1})_+,\quad
\veps_i(b) = \ol{x}_i + (x_{i+1} - \ol{x}_{i+1})_+,\nonumber\\
&&\quad \mbox{for $i=1,2,\cdots,n-1$, and}\nonumber\\
&& \veps_0(b) = l - x_0 - \sum_{i=1}^n (x_i + \ol{x}_i) + 
2(x_1 - \ol{x}_1)_+,\quad
\veps_n(b) = 2 \ol{x}_n +x_0.
\end{eqnarray*}
As before, $(x)_+=\max(x,0)$ and
 $\wt(b)=\sum_{i=0}^n(\vphi_i(b)-\veps_i(b))\La_i$.
In this case the automorphism $\sigma$ is given by
\[
\sigma : (m_0,m_1,\cdots,m_n) \mapsto (m_0,m_1,\cdots,m_n),
\]
for $\la = \sum_{i=0}^n m_i \La_i$.

Choose $b(0) = (0,\cd,0) \in B$, $b(n) = (0,\cd,0,m,x_0,m,0,\cd,0) \in B$, 
where $x_0=\eps(l)$ and $m=\frac12(l-x_0)$. Here $\eps(i)$ is defined in
(\ref{eq:def_eps}).
Then $\vphi(b(i))=l \La_i$, $i=0,n$.
The highest weight vector $u_\la$ with $\la = l \La_i$, $i=0,n$
corresponds to the ground
state path $\pbar = \cd \ot \bbar_k \ot \cd \ot \bbar_1$
where $\bbar_k = b(i)$, $i=0,n$ for all $k$.
Furthermore $\la_k = l \La_i$, $i=0,n$ for all $k$.

\subsection{Description of Demazure crystals}
For $\la = l \La_0,l \ge 1$, set $d=2n$ and choose the sequence $\{i_a^{(j)} |
j \ge 1,1 \le a \le 2n \}$ defined as follows:
\[
i_a^{(j)} = \left\{\begin{array}{ll}
a-1, & 1 \le a \le n+1,\\
2n+1-a, & n+2 \le a \le 2n.
\end{array}
\right.
\]

\begin{theorem}
For $\la = l \La_0$, with the above choice of $i_a^{(j)}$ and $d$,
$B$ satisfies (II),(III) and (IV'). Furthermore, in this case,
\begin{eqnarray*}
&& B_0^{(j)} = \{ (0,\cd ,0) \} \subseteq B, \quad B_{2n}^{(j)} = B, \quad \mbox{and}
\nonumber\\
&& B_a^{(j)} = \{ (x_1,\cd ,x_a,0,\cd,0) \} \subseteq B, \quad 1 \le a
\le n,
\nonumber\\
&& B_{n+a}^{(j)} = \{ (x_1,\cd ,x_n,x_0,\ol{x}_n,
\cd,\ol{x}_{n-a+1},
0,\cd ,0) \} \subseteq B, \quad 1 \le a
\le n-1.
\nonumber\\
\end{eqnarray*}
Here on each set $x_i$'s and $\ol{x}_i$'s run over non negative integers
satisfying the conditions in (\ref{eq:BofD2}).
Also, $b^{(j)}_a$ are given as follows:
\[
b^{(j)}_{0} = (0,\cd,0),
\]
and for $1\leq a\leq n$, 
\[
b^{(j)}_{a} =
(0,\cd ,0,\displaystyle{\mathop{l}_{\scriptstyle a}},0,\cd ,0),
\quad b^{(j)}_{n+a} =
(0,\cd ,0,\displaystyle{\mathop{l}_{\scriptstyle \ol{n-a+1}}},0,\cd ,0).
\]
\end{theorem}

\Proof
The proof is similar to $A^{(2)}_{2n}$.
\qed

\begin{remark}
For $\la = l \La_n$, it can be seen similarly that Theorem \ref{th:iso} holds
for the sequence $\{i_a^{(j)}\}$ given by
\[
i_a^{(j)} = \left\{\begin{array}{ll}
n-a+1, & 1 \le a \le n+1,\\
a-n-1, & n+2 \le a \le 2n.
\end{array}
\right.
\]
\end{remark}

\section{$C_{n}^{(1)}$ case}

\begin{figure}[ht] \label{fig:C}
\hspace{1cm}
\begin{picture}(300,60)(0,0)
\put(5,20){\circle{10}}
\put(10,18){\line(1,0){36}}
\put(10,22){\line(1,0){36}}
\put(50,20){\line(-2,1){10}}
\put(50,20){\line(-2,-1){10}}
\put(55,20){\circle{10}}
\put(60,20){\line(1,0){40}}
\multiput(110,20)(6,0){10}{\circle*{1}}
\put(173,20){\line(1,0){40}}
\put(218,20){\circle{10}}
\put(227,18){\line(1,0){36}}
\put(227,22){\line(1,0){36}}
\put(223,20){\line(2,1){10}}
\put(223,20){\line(2,-1){10}}
\put(268,20){\circle{10}}
\put(0,3){\makebox(10,10){\bf 1}}
\put(50,3){\makebox(10,10){\bf 1}}
\put(213,3){\makebox(10,10){\bf 1}}
\put(263,3){\makebox(10,10){\bf 1}}
\end{picture}
\caption{Dynkin diagram $C^{(1)}_n\; (n \ge 2)$}
\end{figure}
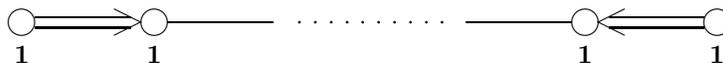

In this section we give explicit descriptions of the $C_{n}^{(1)}$
Demazure crystals $\B_{w}(l\La_i)$, $i=0,n$ for a suitably chosen 
linearly ordered chain of Weyl
group elements $w \in \{ w^{(k)} | k \ge 0 \}$ using the perfect crystal
$B$ of level $l$ (see \cite{KKM}) which is isomorphic to
$B(0) \oplus B(2 \ol{\La}_1) \oplus 
\cdots \oplus B(2 l \ol{\La}_1)$ as crystals for 
$U_q(C_n)$.
It is shown that the mixing index $\kappa = 1$ in these cases.
The Dynkin diagram $C^{(1)}_n\; (n \ge 2)$
is shown in Figure \ref{fig:C}.
The labels are the levels of the fundamental weights corresponding to the 
vertices.

\subsection{Description of the perfect crystal}
For $l \ge 1$, we briefly recall the perfect crystal $B$
(see \cite{KKM}). As a set,
\begin{equation} \label{eq:BofC}
B = \left\{  (x_1,\cdots,x_n,\ol{x}_n,\cdots,\ol{x}_1) \in
\Z^{2n} \Biggm| 
\begin{array}{l}
x_i,\ol{x}_i \ge 0, \\
2l \ge \sum_{i=1}^n
(x_i + \ol{x}_i) \in 2 \Z 
\end{array}
\right\}.
\end{equation}
The action of $\ft{i}$ on $B$ is defined as follows:
For $b=(x_1,\cdots,x_n,\ol{x}_n,
\cdots,\ol{x}_1) \in B$,
\begin{eqnarray*}
&& \ft{0} b = \left\{ 
\begin{array}{ll}
(x_1 +2,x_2,\cdots,\ol{x}_2,\ol{x}_1) & 
\mbox{if $x_1 \ge \ol{x}_1$},\\
(x_1 +1,x_2,\cdots,\ol{x}_2,\ol{x}_1-1) & 
\mbox{if $x_1 = \ol{x}_1-1$},\\
(x_1,x_2,\cdots,\ol{x}_2,\ol{x}_1-2) & 
\mbox{if $x_1 \le \ol{x}_1-2$},
\end{array}
\right.
\nonumber \\
&& \ft{i} b = \left\{ 
\begin{array}{ll}
(x_1,\cdots,x_i-1,x_{i+1}+1,\cdots,\ol{x}_1) & 
\mbox{if $x_{i+1} \ge \ol{x}_{i+1}$},\\
(x_1,\cdots,\ol{x}_{i+1}-1,\ol{x}_i+1,\cdots,\ol{x}_1) & 
\mbox{if $x_{i+1} < \ol{x}_{i+1}$},
\end{array}
\right.
\nonumber \\
&&\quad \mbox{for $i=1,2,\cdots,n-1$, and}\nonumber \\
&& \ft{n} b =(x_1,\cdots,x_n-1,\ol{x}_n+1,\cdots,\ol{x}_1).
\end{eqnarray*}
The action of $\et{i}$ is given by $\et{i}b=b'$ if and only if
$\ft{i}b'=b$. $\vphi_{i}$ and $\veps_{i}$ are given by
\begin{eqnarray*}
&& \vphi_0(b) = l - \frac12 \sum_{i=1}^n (x_i + \ol{x}_i) + 
(\ol{x}_1 - x_1)_+,\quad
\vphi_n(b) = x_n ,\nonumber\\
&& \vphi_i(b) = x_i + (\ol{x}_{i+1} - x_{i+1})_+,\quad
\veps_i(b) = \ol{x}_i + (x_{i+1} - \ol{x}_{i+1})_+,\nonumber\\
&&\quad \mbox{for $i=1,2,\cdots,n-1$, and}\nonumber\\
&& \veps_0(b) = l - \frac12 \sum_{i=1}^n (x_i + \ol{x}_i) + 
(x_1 - \ol{x}_1)_+,\quad
\veps_n(b) = \ol{x}_n .
\end{eqnarray*}
Recall $(x)_+=\max(x,0)$ and $\wt(b)=\sum_{i=0}^n(\vphi_i(b)-\veps_i(b))\La_i$.
In this case the automorphism $\sigma$ is given by
\[
\sigma : (m_0,m_1,\cdots,m_n) \mapsto (m_0,m_1,\cdots,m_n),
\]
for $\la = \sum_{i=0}^n m_i \La_i$.

Choose $b(0) = (0,\cdots,0) \in B$ and 
$\displaystyle
b(n) = (0,\cdots,\mathop{{l}}_n,\mathop{{l}}_{\ol{n}},\cd,0) \in B$.
Then $\vphi(b(i))=l \La_i$, $i=0,n$.
The highest weight vector $u_\la$ with $\la = l \La_i$, $i=0,n$ corresponds to 
the ground state path $\pbar = \cd \ot \bbar_k \ot \cd \ot \bbar_1$,
where $\bbar_k = b(i)$, $i=0,n$ for all $k$.
Furthermore, for all $k$, $\la_k = l \La_i$, $i=0,n$.

\subsection{Description of Demazure crystals}
For $\la = l \La_0,l \ge 1$, set $d=2n$ and choose the sequence $\{i_a^{(j)} |
j \ge 1,1 \le a \le 2n \}$ defined as follows:
\[
i_a^{(j)} = \left\{\begin{array}{ll}
a-1, & 1 \le a \le n+1,\\
2n+1-a, & n+2 \le a \le 2n.
\end{array}
\right.
\]

\begin{theorem}
For $\la = l \La_0$, with the above choice of $i_a^{(j)}$ and $d$, 
$B$ satisfies (II),(III) and (IV').
Furthermore, in this case,
\begin{eqnarray*}
&& B_0^{(j)} = \{ (0,\cd ,0) \} \subseteq B, 
\quad B_{2n}^{(j)} = B, \quad \mbox{and}
\nonumber\\
&& B_a^{(j)} = \{ (x_1,\cd ,x_a,0,\cd,0) \} \subseteq B, \quad 1 \le a
\le n,
\nonumber\\
&& B_{n+a}^{(j)} = \{ (x_1,\cd ,x_n,\ol{x}_n,\cd,\ol{x}_{n-a+1},
0,\cd ,0) \} \subseteq B, \quad 1 \le a
\le n-1.
\nonumber\\
&&
\end{eqnarray*}
Here on each set $x_i$'s and $\ol{x}_i$'s run over non negative integers
satisfying the conditions in (\ref{eq:BofC}).
Also, $b^{(j)}_a$ are given as follows:
\[
b^{(j)}_{0} = (0,\cd,0),
\]
and for $1\leq a\leq n$, 
\[
b^{(j)}_{a} =
(0,\cd ,0,\displaystyle{\mathop{2l}_{\scriptstyle a}},0,\cd ,0),
\quad b^{(j)}_{n+a} =
(0,\cd ,0,\displaystyle{\mathop{2l}_{\scriptstyle \ol{n-a+1}}},0,\cd ,0).
\]
\end{theorem}

\Proof
The proof is similar to $A^{(2)}_{2n}$.
\qed

\begin{remark}
For $\la = l \La_n$, it can be seen similarly that Theorem \ref{th:iso} holds
for the sequence $\{i_a^{(j)}\}$ given by,
\[
i_a^{(j)} = \left\{\begin{array}{ll}
n-a+1, & 1 \le a \le n+1,\\
a-n-1, & n+2 \le a \le 2n.
\end{array}
\right.
\]
\end{remark}

\section{Discussion}

In this section, we discuss three topics: An interesting example from
the $C^{(1)}_n$ case, $\kappa=2$ conjecture, and the classical invariance
of affine Demazure modules. For this purpose, we briefly review the 
situation of the mixing index $\kappa=2$. Let us recall the definition
of $B^{(j)}_a$ ($j\ge1,0\le a\le d$) in section 1.3. We also define
a subset $B^{(j+1,j)}_a$ ($j\ge1,0\le a\le d$) of $B\ot B$ by 
\begin{eqnarray*}
B^{(j+1,j)}_0&=&B^{(j+1)}_0\ot B^{(j)}_d,\\
B_a^{(j+1,j)}&=&\bigcup_{n\ge0}\ft{i_a^{(j+1)}}^n B_{a-1}^{(j+1,j)}
\setminus\{0\}\quad(a=1,\cdots,d).
\end{eqnarray*}
Replace the condition (II) by 
\begin{description}
\item[(II')  ]
For any $j\ge1$,
$B_d^{(j+1,j)}=B_d^{(j+1)}\ot B$.
\end{description}
Note that $B_d^{(j)}\neq B$ in general. Then the corresponding theorem turns
out to be 

\begin{theorem}[{\rm \cite{KMOU}}] 
Under the assumptions (I,II',III,IV), we have
\[
\B_{w^{(k)}}(\la)\simeq
\left\{
\begin{array}{ll}
u_{\la_1}\ot B^{(1)}_a&\mbox{ if }j=1,\\
u_{\la_j}\ot B^{(j,j-1)}_a\ot B^{\ot(j-2)}&\mbox{ if }j\ge2.
\end{array}\right.
\]
\end{theorem}

We move to the first topic. Let us consider the case of $\geh=C^{(1)}_n$.
We set $\la=l\La_i$ ($i\neq0,n$), and try to find a sequence $\{w^{(k)}\}$
for the perfect crystal $B$ in section 8 satisfying $B^{(j)}_d=B$, i.e.,
$\kappa=1$. Take $n=2,l=1,i=1$. Even in this simplest case, it is 
impossible to find such a sequence, but we can easily find a sequence 
satisfying $B^{(j+1,j)}_d=B^{(j+1)}_d\ot B$, i.e., $\kappa=2$. Checking
a few more examples leads to a conjecture that this should be true for any 
$n,l$ and $i$ ($\neq0,n$). Now let us choose another perfect crystal $B'$, 
which is isomorphic to $B(l\ol{\La}_n)$ as $U_q(C_n)$ crystal 
(see \cite{KMN2}).
Then it appears in this case that for $\la=l\La_i$ ($0\le i\le n$) we 
can find a sequence with $\kappa=1$. This shows that the mixing index may 
depend on the choice of the perfect crystals.

We explain what we call $\kappa=2$ conjecture. This originates from a 
natural question how the mixing index changes when we change the highest
weight $\la$ from $l\La_0$ keeping the sequence $\{w^{(k)}\}$ unchanged.
Several simple examples and computations at the character level lead us 
to the following conjecture.

\begin{conjecture}
Let $B$ be a perfect crystal. If the sequence of Weyl group elements
$\{w^{(k)}\}$ satisfies the conditions (II), (III) and (IV) for $\la=l\La_0$,
then the conditions (II') and (III) are also satisfied for any level 
$l$ dominant integral weight $\la$. 
\end{conjecture}
This conjecture has been proved for the case of $\geh=A^{(1)}_n$,
$B=B(l\ol{\La}_1)$ (as $U_q(A_n)$ crystal) in \cite{KMOU}. It is easy
to prove for the cases treated in the remarks of each section.
(Note that in the remarks $\{w^{(k)}\}$ is changed according to the 
change in the highest weights, whereas in the conjecture $\{w^{(k)}\}$
is the same as in the case of $l\La_0$.)

In \cite{KMOTU1}, we have shown the classical
invariance of the Demazure module $V_{w^{(Ld)}}(l\La_0)$ for the case
of $\geh=A^{(1)}_n$. In a similar manner, we can also show the classical
invariance for the cases treated in this paper. Let $V$ be the finite
dimensional $U'_q(\geh)$-module corresponding to the perfect crystal $B$.
Then we have the following isomorphism.
\[
V_{w^{(Ld)}}(l\La_0)\simeq V^{\ot L}
\mbox{ as $U_q(\geh_{I\setminus\{i_L\}})$-modules}.
\]
Here $\geh_{I\setminus\{i\}}$ is the finite dimensional simple Lie
algebra corresponding to the Dynkin diagram obtained by removing the 
vertex $i$, and $i_L$ is determined from $\sigma^L(\La_0)=\La_{i_L}$.
Note that $U_q(\geh_{I\setminus\{i_L\}})$ can be viewed as a subalgebra
of $U_q(\geh)$ in the canonical way.

We have seen in this paper that a perfect crystal $B$ is intimately 
related to a sequence of Weyl group elements $\{w^{(k)}\}$, and the
Demazure crystal $\B_{w^{(k)}}(l\La_0)$ is well described in terms of
paths. Denote by $B^{m,l}$, if it exists, the level $l$ perfect
crystal which is isomorphic to $B(l\gamma_m\ol{\La}_m)\oplus\cd$
as $U_q(\geh_{I\setminus\{0\}})$ crystal.
($\gamma_m$ is some positive integer. In all known cases, $\gamma_m
=\max(1,a_m/a^\vee_m)$, where $a_i$ (resp. $a^\vee_i$) is defined by
$\delta=\sum_{i\in I}a_i\alpha_i$ 
(resp. $c=\sum_{i\in I}a^\vee_i\alpha^\vee_i$).) 
Then it appears that the sequence
$\{w^{(k)}\}$ does not depend on $l$. However, as seen in the $A^{(1)}_n$
case, it is considered that the sequence does depend on $m$. This 
motivates the search for new perfect crystals. We also think it 
important to investigate Demazure characters and clarify their relation
to Kostka type polynomials and $q$-analogues of products of the classical
characters. In \cite{KMOTU1}, we have picked the $A^{(1)}_n$ case,
and discussed the relation between the Demazure character and the product
of Schur functions or the Kostka-Foulkes polynomial. We hope the classical
invariance discussed here shed new light on the investigation of these
polynomials. We give a unified treatment of level one Demazure characters in 
terms of one dimensional configuration sums in another publication 
\cite{KMOTU2}.

\bigskip
\noindent{\em Acknowledgement.}\quad
A.K. thanks T. Nakanishi and A. Varchenko for their warm hospitality at
University of North Carolina. K.C.M. and M.O. thank O. Foda for his 
hospitality at the University of Melbourne, and for collaboration in
their earlier work \cite{FMO}. K.C.M is partially supported by NSA/MSP 
Grant No. 96-1-0013.

\end{document}